\documentclass{amsart}
\usepackage{amssymb}
\usepackage{amsfonts}
\usepackage{amstext}
\usepackage{algorithmic}
\usepackage{algorithm}
\usepackage{graphicx}
\usepackage{epstopdf}
\usepackage[all]{xy}

\allowdisplaybreaks

\numberwithin{equation}{section}
\newtheorem{theorem}{Theorem}[section]

\newtheorem{corollary}[theorem]{Corollary}

\theoremstyle{definition}

\newtheorem{example}[theorem]{Example}

\theoremstyle{remark}
\newtheorem{remark}[theorem]{\bf{Remark}}

\newcommand{\Ccal}{{\mathcal{C}}}

\newcommand{\Mcal}{{\mathcal{M}}}

\newcommand{\Rcal}{{\mathcal{R}}}

\newcommand{\tens}{\otimes}
\newcommand{\id}{{\rm id}}
\newcommand{\bo}{{}^{(1)}}
\newcommand{\bt}{{}^{(2)}}
\renewcommand{\o}{{}_{(1)}}
\renewcommand{\t}{{}_{(2)}}

\newcommand{\la}{{\triangleright}}

\newcommand{\bos}{\rtimes\kern -5.9 pt\cdot}
\newcommand{\crosscop}{\blacktriangleright\kern -4 pt\blacktriangleleft}

\begin{document}

\title{Transmutation and Bosonisation of Quasi-Hopf Algebras}

\author{Jennifer Klim}
\address{Queen Mary, University of London\\
School of Mathematics, Mile End Rd, London E1 4NS, UK}

\date{\today}

\begin{abstract}

Let $H$ be a quasitriangular quasi-Hopf algebra, we construct a braided group $\underline{H}$ in the quasiassociative category of left $H$-modules. Conversely, given any braided group $B$ in this category, we construct a quasi-Hopf algebra $B\bos H$ in the category of vector spaces. We generalise the transmutation and bosonisation theory of \cite{Majid94} to the quasi case. As examples, we bosonise the octonion algebra to an asoociative one, obtain the twisted quantum double $D^{\phi}(G)$ of a finite group as a bosonisation, and obtain its transmutation. Finally, we show that $\underline{H}\bos H \cong H _\Rcal\crosscop H$ as quasi-Hopf algebras.
\end{abstract}

\maketitle

\section{Introduction}

If $H$ is a quasitriangular Hopf algebra, it is known that there exists a Hopf algebra $\underline{H}$ in the category $_H\Mcal$ of left $H$-modules, a construction known as `transmutation' \cite{Majid93}. Following Majid, we refer to Hopf algebras in braided categories as `braided groups'. Conversely, given a braided group $B$ in the category $_H\Mcal$, there exists a Hopf algebra $B\bos H$ in the category of vector spaces \cite{Majid94}, a construction known as `bosonisation'. We recall the required theory in section 2.

In sections 3 and 4 we generalise these results to $H$ a quasitriangular quasi-Hopf algebra \cite{Drinfeld90}. In this case the associativity constraint in the category $_H\Mcal$ is no longer trivial, it now depends on the associator of the quasi-Hopf algebra. Nevertheless, we show that the theory of \cite{Majid94} follows through. One has a transmutation $\underline{H}$ as a braided group in $_H\Mcal$. In \cite{Bulacu00} it was shown that for any algebra $B$ in $_H\Mcal$ there is an associative algebra $B\bos H$. We extend this to $B$ a braided group in $_H\Mcal$ and obtain a quasi-Hopf algebra $B\bos H$. One also has, for example a one to one correspondence between braided $B$-modules in $_H\Mcal$ and left $B\bos H$-modules in the category of vector spaces. We consider the examples of the twisted quantum double $D^{\phi}(G)$ introduced in \cite{Dijkgraaf90}, and the octonions in the form \cite{Majid00}.

It is known for quasitriangular Hopf algebras that there exists an isomorphism between $\underline{H}\bos H$ and the twisted tensor product $H _{\Rcal}\crosscop H$; in section 5 we prove that when $H$ is a quasitriangular quasi-Hopf algebra, there is a quasi-Hopf algebra isomorphism $\chi:\underline{H}\bos H\to H _{\Rcal}\crosscop H$.

\section{Preliminaries}

\subsection{Quasi-Hopf Algebras}

Let $k$ be a commutative field. A \textit{quasi-bialgebra}, \cite{Drinfeld90}, over $k$ is $(H, \Delta, \varepsilon, \phi)$ where $H$ is a unital associative algebra over $k$, $\Delta: H \rightarrow H \otimes H$ is an algebra homomorphism such that

\begin{equation*}(\id \otimes \Delta)\Delta = \phi(\Delta \otimes \id)\Delta\phi^{-1},\end{equation*}

and the axiom for the counit $\varepsilon$, an algebra homomorphism, are as usual. The element $\phi \in H \otimes H \otimes H$, called the \emph{Drinfeld associator}, or \emph{associator}, that controls the noncoassociativity is invertible, and is required to be a counital 3-cocyle, i.e.

\begin{equation*}(1 \otimes \phi)(\id \otimes \Delta \otimes \id)\phi(\phi \otimes 1) = (\id \otimes \id \otimes \Delta)\phi(\Delta \otimes \id \otimes \id)\phi,\end{equation*}

and $(\id \otimes \varepsilon \otimes \id)(\phi) = 1 \otimes 1 \otimes 1$. $H$ is a \textit{quasi-Hopf algebra} if there exists a convolution invertible algebra anti-homomorphism $S: H \rightarrow H$, called the antipode, together with elements $\alpha ,\beta \in H$ such that,

\begin{equation*} S(h\o) \alpha h\t = \varepsilon(h) \alpha, \end{equation*}
\begin{equation*} h\o \beta S(h\t) = \varepsilon(h) \beta, \end{equation*}
\begin{equation*} X^1 \beta S(X^2) \alpha X^3 = 1, \end{equation*}
\begin{equation*} S(x^1) \alpha x^2 \beta S(x^3) = 1, \end{equation*}

for all $h \in H$, where $\phi =  X^1 \otimes X^2 \otimes X^3$ is written in capital letters, and $\phi^{-1} =  x^1 \otimes x^2 \otimes x^3$ is written in lower case letters. For brevity, the sum notation for the coproduct and the Drinfeld associator has been suppressed. The antipode is uniquely determined up to a transformation $\alpha \mapsto U\alpha$, $\beta \mapsto \beta U^{-1}$, $S(h) \mapsto US(h)U^{-1}$, for any invertible $U\in H$. Following from this, we can, without loss of generality, assume $\varepsilon(\alpha) = \varepsilon(\beta) = 1$.

For Hopf algebras it is known that the antipode is a coalgebra anti-homomorphism; in the case of quasi-Hopf algebras this is true only up to a twist, i.e. there exists $f\in H\tens H$ such that

\[ f\Delta(S(h))f^{-1}=S(\Delta^{op}(h))\]

for all $h\in H$.

Following \cite{Drinfeld90}, define $\gamma,\delta\in H\tens H$ by

\[ \gamma = S(A^2)\alpha A^3 \tens S(A^1)\alpha A^4 \]
\[ \delta = B^1\beta S(B^4)\tens B^2\beta S(B^3) \]

where

\[ A^1\tens A^2\tens A^3\tens A^4 = (\phi\tens 1)(\Delta\tens\id\tens\id)(\phi^{-1}) \]
\[ B^1\tens B^2\tens B^3\tens B^4 = (\Delta\tens\id\tens\id)(\phi)(\phi^{-1}\tens 1) \]

Denote $f^{-1}$ by $g$, then $f,g$ are given by the formulae

\[ f = (S\tens S)(\Delta^{op}(x^1))\gamma\Delta(x^2\beta S(x^3)) \]
\[ g = \Delta(S(x^1)\alpha x^2)\delta (S\tens S)(\Delta^{op}(x^3)) \]

Further, $f$ satisfies $f\Delta(\alpha) = \gamma$, $\Delta(\beta)g=\delta$, and we note

\[ \Delta(X^1)\delta(S\tens S)(\Delta^{op}(X^2))\gamma \Delta(X^3) = 1 \]
\[ (S\tens S)(\Delta^{op}(x^1))\gamma\Delta(x^2)\delta(S\tens S)(\Delta^{op}(x^3)) =1 \]

It is useful to define elements $q = q^1 \otimes q^2 = \sum X^1 \otimes S^{-1}(\alpha X^3) X^2$ and $p = p^1 \otimes p^2 = \sum x^1 \otimes x^2 \beta S(x^3)$ in $H \otimes H$. Then, for all $h \in H$,

\begin{equation*} \Delta(h\o) p (1 \otimes S(h\t)) = p (h \otimes 1), \end{equation*}
\begin{equation*} (1 \otimes S^{-1}(h\t)) q \Delta(h\o) = (h \otimes 1) q, \end{equation*}
\begin{equation*} \Delta(q^1) p (1 \otimes S(q^2)) = 1 \otimes 1, \end{equation*}
\begin{equation*} (1 \otimes S^{-1}(p^2)) q \Delta(p^1) = 1 \otimes 1. \end{equation*}

The quasi-Hopf algebra $(H, \Delta, \varepsilon, S, \alpha, \beta, \phi)$ is \textit{quasitriangular} \cite{Drinfeld90} if there is an invertible element $R \in H \otimes H$ such that,

\begin{equation*} (\Delta \otimes \id)(R) = \phi_{312}R_{13}\phi^{-1}_{132}R_{23}\phi, \end{equation*}
\begin{equation*} (\id \otimes \Delta)(R) = \phi^{-1}_{231}R_{13}\phi_{213}R_{12}\phi^{-1}, \end{equation*}
\begin{equation*} \Delta^{op}(h) = R\Delta(h)R^{-1}, \end{equation*}

for all $h \in H$. Writing $\phi = \sum X^1 \otimes X^2 \otimes X^3$, then $\phi_{ijk} \in H \otimes H \otimes H$ has $X^1$ in the i-th position, $X^2$ in the j-th position and $X^3$ in the k-th position, for example, $\phi_{312} = X^2 \otimes X^3 \otimes X^1$. Similarly for $\phi = \sum x^1 \otimes x^2 \otimes x^3$. The inverse \cite{Bulacu03} is given by

\begin{equation*} R^{-1}=X^1\beta S(Y^2R\bo x^1X^2)\alpha Y^3x^3X^3\t\tens Y^1R\bt x^2X^3\o \end{equation*}

As for quasitriangular Hopf algebras, $(\varepsilon \otimes \id)(R) = (\id \otimes \varepsilon)(R) = 1 \otimes 1$. Further, the above relations imply the quasi-Yang-Baxter equation

\[ R_{12}\phi_{312}R_{13}\phi^{-1}_{132}R_{23}\phi = \phi_{321}R_{23}\phi^{-1}_{231}R_{13}\phi_{213}R_{12}. \]

\subsection{Monoidal Categories}

A \textit{monoidal category} is $(\mathcal{C}, \otimes, I, \Phi, l, r)$, where $\mathcal{C}$ is a category, $\otimes: \mathcal{C} \times \mathcal{C} \rightarrow \mathcal{C}$ is a functor called the tensor product, and $I$ is a fixed unit object. Further, $\mathcal{C}$ is equipped with a natural transformation, called the \emph{associativity constraint}, consisting of functorial isomorphisms $\Phi_{U,V,W}:(U \otimes V) \otimes W \rightarrow U \otimes (V \otimes W)$ for all $U,V,W \in \mathcal{C}$, obeying the well-known pentagon condition. Finally, $\mathcal{C}$ is equipped with  natural transformations $l,r$, consisting of functorial isomorphisms $l_V: V \rightarrow I \otimes V $, $r_V: V \rightarrow V \otimes I$ for each $V \in \mathcal{C}$, obeying the triangle condition. Here $l$ and $r$ are called the \emph{left \emph{and} right unit constraints} respectively.

An object $V$ in a monoidal category $\mathcal{C}$ is \textit{rigid} if there exists an object $V^*$ and morphisms $ev_V:V^* \otimes V \rightarrow I$ and $coev_V:I \rightarrow V \otimes V^*$ in $\mathcal{C}$, such that

\[ r^{-1}_V (\id_V \otimes ev_V)\Phi_{V,V^*,V}(coev_V \otimes \id_V) l_V = \id_V, \]
\[ l^{-1}_{V^*} (ev_V \otimes \id_{V^*})\Phi^{-1}_{V^*,V,V^*}(\id_{V^*} \otimes coev_V)r_{V^*} = \id_{V^*}. \]

The monoidal category $\mathcal{C}$ is called \emph{rigid} if every object in $\mathcal{C}$ is rigid. A \emph{braided} category \cite{Joyal93} is a monoidal category  $(\mathcal{C}, \otimes)$ equipped with a natural transformation consisting of functorial isomorphisms $\Psi_{V,W}:V \otimes W \rightarrow W \otimes V$ for all $V,W \in \mathcal{C}$, called a \emph{braiding}, obeying the well-known hexagon conditions. We will use the notation of \cite{Majid95}.

\begin{example} \cite{Drinfeld90} Let $H$ be a unital algebra, then the category $_{H}\mathcal{M}$ of left $H$-modules consists of objects, the vector spaces $V$ on which $H$ acts, and morphisms, the linear maps $f$ which commute with the action of $H$, i.e. $f(h \la v) = h \la f(v)$ for all $v \in V$ and $h \in H$. If $H$ is a quasi-bialgebra, then $\otimes$, defined by $h \la (v \otimes w) = h\o \la v \otimes h\t \la w$, and

\[ \Phi_{U,V,W}((u \otimes v) \otimes w) = X^1 \la u \otimes (X^2 \la v \otimes X^3 \la w) \]

for all $u \in U, v \in V, w \in W$ where $U,V,W \in {_H}\mathcal{M}$, makes $_{H}\mathcal{M}$ into a monoidal category. If $H$ is a quasi-triangular quasi-Hopf algebra, then $_{H}\mathcal{M}$ is a braided monoidal category with the braiding defined by

\[ \Psi_{U,V}(u \otimes v) = R\bt \la v \otimes R\bo \la u \]

for all $u \in U, v \in V$. Finally, this category is rigid with $(h \la f)(v) = f(S(h) \la v)$ for all $v \in V$, $f \in V^*$ and $h \in H$, and

\[ ev(f \otimes v) = f(\alpha \la v) \]
\[ coev = \sum_a \beta \la e_a \otimes f^a \]

where $\{e_a\}$ is a basis for $V$ and $\{f^a\}$ a dual basis. We refer to \cite{Majid95} for details.

\end{example}

\subsection{Hopf Algebras in Braided Categories}

An algebra in a monoidal category $\mathcal{C}$ is an object $B$ of $\mathcal{C}$ equipped with a multiplication morphism $B \otimes B \rightarrow B$ and a unit morphism $\underline{1} \rightarrow B$, obeying the usual associativity and unit axioms, but now as morphisms in $\mathcal{C}$, and where $B \otimes B$ is the tensor product in the catgeory. A bialgebra in a braided category is an algebra $B$ in the category equipped with algebra morphisms $\underline{\Delta}: B \rightarrow B \otimes B$ and $\underline{\varepsilon}: B \rightarrow \underline{1}$ in $\mathcal{C}$ which form a coalgebra in the category. Further, if there is a morphism $\underline{S}:B \rightarrow B$ in $\mathcal{C}$ obeying the usual antipode axioms, then $B$ is a Hopf algebra in the braided category $\mathcal{C}$. The Hopf algebra $B$ is called a \emph{braided Hopf algebra} or \emph{braided group}, \cite{Majid93a}.

Following \cite{Majid95}, we consider monoidal categories $\mathcal{C}$ and $\mathcal{D}$ with $\mathcal{D}$ braided, and functors $F, V\tens F:\mathcal{C} \rightarrow \mathcal{D}$. Suppose there is an object $B \in \mathcal{D}$ such that for all $V \in \mathcal{D}$, $Mor(V,B) \cong Nat(V \otimes F,F)$ by functorial isomorphisms $\theta_V$. Let

\[ \alpha = \{ \alpha_M : B \otimes F(M) \rightarrow F(M) | M \in \mathcal{C} \} \]

be the natural transformation corresponding to the identity morphism $\id_B$ in $Mor(V,B)$. Then, using $\alpha$, and the braiding we get induced maps

\[ \theta^n_V : Mor(V,B^{\otimes n}) \rightarrow Nat(V \otimes F^n , F^n) \]

and we assume these are bijections. This is called the \textit{representability assumption for modules}. Then, using these bijections, we can define a multiplication, a unit, a coproduct, a counit and an antipode for $B$.

For example, note that $\alpha_M  (\id \otimes \alpha_M)\Phi_{B,B,F(M)}: (B \otimes B) \otimes F(M) \rightarrow F(M)$ is a natural transformation in $Nat(B \otimes B \otimes F,F)$, and hence corresponds to a unique map $B \otimes B \rightarrow B$ under $\theta^{-1}_{B \otimes B}$, which must be the multiplication on $B$. We will require the following theorem,

\begin{theorem} \cite{Majid95} Let $\mathcal{C}$ and $\mathcal{D}$ be monoidal categories with $\mathcal{D}$ braided, and $F: \mathcal{C} \rightarrow \mathcal{D}$ be a monoidal functor satisfying the representability assumption for modules. Then $B$, as above, is a bialgebra in $\mathcal{D}$. If $\mathcal{D}$ is rigid, then $B$ is a Hopf-algebra in $\mathcal{D}$.
\end{theorem}

This theorem with $\Ccal = _H\Mcal$, for $H$ a quasitriangular Hopf algebra, and $F=\id$ is used to reconstruct a braided group, $\underline{H}$ \cite{Majid94}. Taking the monoidal category of $B$-modules in the braided category of $H$-modules and the forgetful functor to $Vec$, and reconstructing, we obtain an ordinary Hopf algebra, which is the categorical theory of bosonisation. We now do the same when $H$ is a quasitriangular quasi-Hopf algebra.

\section{Transmutation of Quasi-Hopf Algebras}

Let $H$ be a quasitriangular quasi-Hopf algebra, $B_L$ be the vector space $H$ with the left regular action, and let $B$ be the vector space $H$ viewed as an object of $_H \mathcal{M}$ via the left adjoint action $h \la g = h\o gS(h\t)$ for all $h,g \in H$. In the notation of the above theorem, we consider the case when $\mathcal{C} = \mathcal{D} = {_H} \mathcal{M}$, and $F=\id$.

First we define $\theta_V:Mor(V,B)\rightarrow Nat(V \otimes \id, \id)$ for $V \in {_H}\mathcal{M}$ as follows. Given $\psi \in Mor(V,B)$, we define $\xi \in Nat(V \otimes \id,\id)$ by

\[ \xi_M(v \otimes m) = \theta_V(\psi)_M(v \otimes m) = q^1 \psi(v) S(q^2) \la m, \]

where $\la$ is the action of $H$ on $M$ as an object in the category $_H \mathcal{M}$. We have to check that each $\xi_M:V \otimes M \rightarrow M$ is a morphism in the category if $\psi$ is.

\begin{eqnarray*}
\xi_M(h \la (v \otimes m)) & = & \xi_M(h\o \la v \otimes h\t \la m) \\
 & = & q^1 \psi(h\o \la v) S(q^2) h\t \la m \\
 & = & q^1 (h\o \la \psi(v)) S(q^2) h\t \la m \\
 & = & q^1 h_{(1,1)} \psi(v) S(S^{-1}(h\t)q^2h_{(1,2)}) \la m \\
 & = & h q^1 \psi(v) S(q^2) \la m \\
 & = & h \la (q^1 \psi(v) S(q^2) \la m) \\
 & = & h \la (\xi_M(v \otimes m))
\end{eqnarray*}

Conversely, we define $\theta^{-1}_V:Nat(V \otimes \id, \id) \rightarrow Mor(V,B)$ for $V \in {_H}\mathcal{M}$ as follows. Given $\xi \in Nat(V \otimes \id, \id)$, we define $\psi \in Mor(V,B)$ by

\[ \psi(v) = \theta^{-1}_V(\xi)(v) = \xi_{B_L}(p^1 \la v \otimes p^2), \]

for all $v \in V$. Now,

\begin{eqnarray*}
h \la \psi(v) & = & h \la \xi_{B_L}(p^1 \la v \otimes p^2) \\
 & = & h\o \xi_{B_L}(p^1 \la v \otimes p^2) S(h\t) \\
 & = & \xi_{B_L} (h\o \la (p^1 \la v \otimes p^2)) S(h\t) \\
 & = & \xi_{B_L} (h_{(1,1)}p^1 \la v \otimes h_{(1,2)} \la p^2) S(h\t) \\
 & = & \xi_{B_L} (h_{(1,1)}p^1 \la v \otimes h_{(1,2)}p^2) S(h\t) \\
 & = & \xi_{B_L} (h_{(1,1)}p^1 \la v \otimes h_{(1,2)}p^2 S(h\t)) \\
 & = & \xi_{B_L} (p^1 h \la v \otimes p^2) \\
 & = & \xi_{B_L} (p^1 \la (h \la v) \otimes p^2) \\
 & = & \psi(h \la v).
\end{eqnarray*}

It is straightfoward to check these two processes are mutually inverse. The natural transformation corresponding to the identity morphism on $B$ is $\alpha = \{\alpha_M | M \in {_H}\mathcal{M}\}$, where each $\alpha_M :B \otimes M \rightarrow M$ is given by

\[ \alpha_{M}(b \otimes m) = \theta_{B}(\id_B)_M(b \otimes n) = q^1bS(q^2) \la m \]

\begin{theorem} Every quasitriangular quasi-Hopf algebra $H$ has a braided group analogue $\underline{H}$ in $_H\Mcal$ and is given by

\[ \underline{m}(b \otimes b') = q^1(x^1\la b)S(q^2)x^2b'S(x^3) \]
\[ \underline{\eta}(1) = \beta \]
\[ \underline{\Delta}(b) = x^1X^1b\o g^1S(x^2R\bt y^3X^3\t)\tens x^3R\bo\la y^1X^2b\t g^2S(y^2X^3\o) \]
\[ \underline{\varepsilon}(b) = \varepsilon(b) \]
\[ \underline{S}(b) = X^1R\bt x^2\beta S(q^1(X^2R\bo x^1 \la b)S(q^2)X^3x^3) \]

\end{theorem}

\proof
Let $M \in {_H}\mathcal{M}$. We have that $\alpha_M (\id \tens \alpha_M) \Phi_{B,B,M}: (B \tens B) \tens M \rightarrow M$ is a natural transformation in $Nat(B \tens B \tens \id, \id)$, and hence corresponds to a unique map $\underline{m}:B \tens B \rightarrow B$ under $\theta^{-1}_{B \tens B}$. Let $\xi_M = \alpha_M (\id \tens \alpha_M) \Phi_{B,B,M}$, then for all $b,b' \in B$ and $m \in M$,

\begin{eqnarray*}
\xi_M((b \tens b') \tens m) & = & \alpha_M(\id \tens \alpha_M)\Phi_{B,B,M}((b \tens b') \tens m) \\
 & = & \alpha_M(\id \tens \alpha_M)(X^1 \la b \tens (X^2 \la b' \tens X^3 \la m)) \\
 & = & \alpha_M (X^1 \la b \tens q^1(X^2 \la b')S(q^2)X^3 \la m) \\
 & = & Q^1(X^1\la b)S(Q^2)q^1(X^2\la b')S(q^2)X^3 \la m
\end{eqnarray*}

where $Q^1 \tens Q^2$ is another copy of $q = q^1 \tens q^2$. Then,

\begin{eqnarray*}
\underline{m}(b \tens b') & = & \theta^{-1}_{B \tens B}(\xi)(b \tens b') \\
 & = & \xi_{B_L} (p^1 \la (b \tens b') \tens p^2) \\
 & = & \xi_{B_L}((p^1\o \la b \tens p^1\t \la b') \tens p^2) \\
 & = & Q^1(X^1 p^1\o \la b)S(Q^2)q^1(X^2 p^1\t \la b')S(q^2) X^3 \la p^2 \\
 & = & Q^1(X^1 p^1\o \la b)S(Q^2)q^1(X^2 p^1\t \la b')S(q^2) X^3 p^2 \\
 & = & Q^1(X^1 x^1\o \la b)S(Q^2)q^1(X^2 x^1\t \la b')S(q^2) X^3 x^2 \beta S(x^3) \\
 & = & Q^1(x^1 X^1 \la b)S(Q^2)q^1(x^2\o y^1 X^2 \la b')S(q^2)x^2\t y^2X^3\o \\
 &   & \beta S(X^3\t)S(X^3y^3) \\
 & = & Q^1(x^1X^1 \la b)S(Q^2)q^1(x^2\o y^1X^2 \la b')S(q^2)x^2\t y^2\varepsilon(X^3) \\
 &   & \beta S(x^3y^3) \\
 & = & Q^1(x^1 \la b)S(Q^2)q^1(x^2\o y^1 \la b')S(q^2)x^2\t y^2 \beta S(x^3y^3) \\
 & = & Q^1(x^1 \la b)S(Q^2)q^1x^2_{(1,1)}(y^1 \la b')S(S^{-1}(x^2\t q^2 x^2_{(1,2)})y^2 \\
 &   & \beta S(x^3y^3) \\
 & = & Q^1(x^1 \la b)S(Q^2)x^2q^1(y^1 \la b')S(q^2)y^2\beta S(x^3y^3) \\
 & = & X^1(x^1 \la b)S(S^{-1}(\alpha X^3)X^2)x^2q^1(y^1 \la b')S(q^2)y^2 \\
 &   & \beta S(x^3y^3) \\
 & = & X^1(x^1 \la b)S(X^2)\alpha X^3x^2q^1(y^1\la b')S(q^2)y^2 \\
 &   & \beta S(x^3y^3) \\
 & = & X^1x^1\o b S(X^2x^1\t)\alpha X^3x^2q^1(y^1 \la b')S(q^2)y^2 \\
 &   & \beta S(x^3y^3) \\
 & = & X^1x^1\o b S(X^2x^1\t)\alpha X^3x^2q^1(p^1 \la b')S(q^2)p^2S(x^3) \\
 & = & y^1X^1 b S(y^2\o x^1X^2)\alpha y^2\t x^2X^3\o q^1(p^1 \la b')S(q^2) \\
 &   & p^2S(y^3x^3X^3\t) \\
 & = & y^1X^1b S(x^1X^2)S(y^2\o)\alpha y^2\t x^2X^3\o q^1(p^1 \la b')S(q^2) \\
 &   & p^2S(y^3x^3X^3\t) \\
 & = & y^1X^1b S(x^1X^2)\varepsilon(y^2)\alpha x^2X^3\o q^1(p^1 \la b')S(q^2) \\
 &   & p^2S(y^3x^3X^3\t) \\
 & = & X^1b S(x^1X^2)\alpha x^2X^3\o q^1(p^1 \la b')S(q^2)p^2S(x^3X^3\t) \\
 & = & X^1b S(x^1X^2) \alpha x^2X^3\o q^1p^1\o b' S(S^{-1}(p^2)q^2p^1\t) S(x^3X^3\t) \\
 & = & X^1 b S(x^1X^2) \alpha x^2X^3\o b' S(x^3X^3\t)\\
 & = & X^1x^1\o b S(X^2x^1\t)\alpha X^3x^2b'S(x^3)\\
 & = & q^1(x^1\la b)S(q^2)x^2b'S(x^3)
\end{eqnarray*}

So, for all $b,b' \in B$, the multiplication is defined by

\[ \underline{m}(b \tens b') = q^1(x^1\la b)S(q^2)x^2b'S(x^3). \]
 
The antipode is determined by $\underline{S}(b)=\theta^{-1}_B(\xi)(b)$, where 

$\xi_M=r^{-1}_M(M \tens ev_M)\Phi_{M,M^*,M}((M \tens \alpha_{M^*})\tens M)(\Phi_{M,B,M^*} \tens M) ((\Psi_{B,M} \tens M^*)\tens M)(\Phi^{-1}_{B,M,M^*} \tens M)((B \tens coev_M) \tens M) (r_B \tens M)$.

So,

\[ \xi_M(b\tens m)	=	 Q^1X^1R\bt x^2\beta S(q^1(X^2R\bo x^1 \la b)S(q^2)X^3x^3)S(Q^2) \la m \]

hence, 

\begin{eqnarray*}
\underline{S}(b)	&	=	& \theta^{-1}_B(\xi)(b) \\
									&	=	&	\xi_{B_L}(p^1\la b \tens p^2) \\
									&	=	&	Q^1X^1R\bt x^2\beta S(q^1(X^2R\bo x^1p^1 \la b)S(q^2)X^3x^3)S(Q^2) \la p^2\\
									&	=	& Q^1X^1R\bt P^2 S(q^1(X^2R\bo P^1p^1\la b)S(q^2)X^3x^3)S(Q^2)\\
									&	=	&	Q^1X^1R\bt p^1\o\t P^2S(p^1\t)S(q^1(X^2R\o p^1\o\o P^1\la b)S(q^2)X^3)S(Q^2)p^2\\
									&	=	& Q^1X^1p^1\o\o R\bt P^2S(q^1(X^2p^1\o\t R\bo P^1\la b)S(q^2)X^3p^1\t)S(Q^2)p^2\\
									&	=	&	Q^1p^1\o X^1R\bt P^2S(q^1(p^1\t\o X^2 R\bo P^1\la b)S(q^2)p^1\t\t X^3)S(Q^2)p^2\\
									&	=	&	Q^1p^1\o X^1R\bt P^2S(p^1\t q^1(X^2R\bo P^1\la b)S(q^2)X^3)S(Q^2)p^2\\
									&	=	& Q^1p^1\o X^1R\bt P^2S(q^1(X^2R\bo P^1\la b)S(q^2)X^3)S(S^{-1}(p^2)Q^2p^1\o)\\
									&	=	& X^1R\bt p^2S(q^1(X^2R\bo p^1\la b)S(q^2)X^3)
\end{eqnarray*}

So the antipode is defined as

\[ \underline{S}(b) = X^1R\bt p^2S(q^1(X^2R\bo p^1\la b)S(q^2)X^3) \]

The reconstructed $\underline{\Delta}$ is characterised by

\begin{equation*}
\begin{split}
\alpha_{M \tens N}\Phi_{B,M,N}  = & (\alpha_M \tens \alpha_N)\Phi^{-1}_{B,M,B \tens N}(B \tens \Phi_{M,B,N})(B \tens (\Psi_{B,M} \tens N))\\
 & (B \tens \Phi^{-1}_{B,M,N})\Phi_{B,B,M \tens N}(\underline{\Delta} \tens (M \tens N))\Phi_{B,M,N}.
\end{split}
\end{equation*}

Let $b \in B, m \in M, n \in N$, then

\begin{eqnarray*}
\alpha_{M \tens N} \Phi_{B,M,N} ((b \tens m) \tens n) & = & \alpha_{M \tens N} (X^1 \la b \tens (X^2 \la m \tens X^3 \la n)) \\
 & = & q^1(X^1 \la b)S(q^2) \la (X^2 \la m \tens X^3 \la n) \\
 & = & (q^1(X^1 \la b)S(q^2))\o X^2 \la m \tens (q^1(X^1 \la b)S(q^2))\t X^3 \la n
\end{eqnarray*}

and,

\begin{eqnarray*}
\lefteqn{(\alpha_M \tens \alpha_N) \ldots  \Phi_{B,M,N}((b \tens m) \tens n)} \\
 & = & q^1(y^1Y^1\la(X^1 \la b)\underline{\o})S(q^2)y^2Z^1R\bt x^2Y^3\o X^2 \la m \tens \\
 &   &  Q^1(y^3\o Z^2R\bo x^1Y^2\la(X^1 \la b)\underline{\t})S(Q^2)y^3\t Z^3x^3Y^3\t X^3 \la n
\end{eqnarray*}

Since these are equal for all $b \in B, m \in M, n \in N$, we have

\begin{eqnarray*}
\lefteqn{(q^1(X^1 \la b)S(q^2))\o X^2 \tens (q^1(X^1 \la b)S(q^2))\t X^3} \\
 & = & q^1(y^1Y^1\la(X^1 \la b)\underline{\o})S(q^2)y^2Z^1R\bt x^2Y^3\o X^2  \tens \\
 &   & Q^1(y^3\o Z^2R\bo x^1Y^2\la(X^1 \la b)\underline{\t})S(Q^2)y^3\t Z^3x^3Y^3\t X^3
\end{eqnarray*}

Which can be further simplified to
\begin{eqnarray}
\lefteqn{\Delta( q^1bS(q^2) ) = q^1(y^1X^1 \la b\underline{\o})S(q^2)y^2Y^1R\bt x^2 X^3\o } \nonumber \\
  &   & \,\,\,\,\,\,\,\,\,\, \tens Q^1(y^3\o Y^2R\bo x^1X^2 \la b\underline{\t})S(Q^2)y^3\t Y^3x^3 X^3\t \label{comult}
\end{eqnarray}

We can check that $\underline{\Delta}(b) = x^1X^1b\o g^1S(x^2R\bt y^3X^3\t)\tens x^3R\bo\la y^1X^2b\t g^2S(y^2X^3\o)$ satisfies this identity as follows.

\begin{eqnarray*}
\lefteqn{q^1(y^1X^1\la b\underline{\o})S(q^2)y^2Y^1R\bt x^2X^3\o\tens Q^1(y^3\o Y^2R\bo x^1X^2\la b\underline{\t})S(Q^2)y^3\t Y^3x^3X^3\t}\\
	&	=	& q^1(y^1X^1\la w^1A^1b\o g^1S(w^2R'\bt z^3A^3\t))S(q^2)y^2Y^1R\bt x^2X^3\o\\
	&		&	\tens Q^1(y^3\o Y^2R\bo x^1X^2w^3R'\bo\la z^1A^2b\t g^2S(z^2A^3\o))S(Q^2)y^3\t Y^3x^3X^3\t \\
	&	=	&	W^1y^1\o \underline{X^1\o w^1}A^1b\o g^1S(W^2y^1\t \underline{X^1\t w^2}R'\bt z^3A^3\t)\alpha W^3y^2R\bt x^2\underline{X^3\o} \\
	&		& \tens Q^1(y^3\o Y^2R\bo x^1\underline{X^2w^3}R'\bo\la z^1A^2b\t g^2S(z^2A^3\o))S(Q^2)y^3\t Y^3x^3\underline{X^3\t}\\
	&	=	& W^1y^1\o w^1X^1A^1b\o g^1S(W^2y^1\t w^2T^1X^2\o R'\bt z^3A^3\t)\alpha \\
	&		& \,\,\,\,\, W^3y^2\underline{Y^1R\bt x^2w^3\t\o} T^3\o X^3\o \\
	&		& \tens Q^1(y^3\o\underline{Y^2R\bo x^1w^3\o} T^2X^2\t R'\bo\la z^1A^2b\t g^2S(z^2A^3\o))S(Q^2)\\
	&		& \,\,\,\,\, y^3\t\underline{Y^3x^3w^3\t\t} T^3\t X^3\t\\
	&	=	& \underline{W^1y^1\o w^1}X^1A^1b\o g^1S(\underline{W^2y^1\t w^2}T^1X^2\o R'\bt z^3A^3\t)\alpha\\
	&		& \,\,\,\,\, \underline{W^3y^2w^3\o} Y^1R\bt x^2T^3\o X^3\o  \\
	&		& \tens Q^1(\underline{y^3\o w^3\t\o} Y^2R\bo x^1 T^2X^2\t R'\bo\la z^1A^2b\t g^2S(z^2A^3\o))S(Q^2)\\
	&		& \,\,\,\,\, \underline{y^3\t w^3\t\t} Y^3x^3T^3\t X^3\t\\
	&	=	& y^1X^1A^1b\o g^1S(y^2\o w^1T^1X^2\o R'\bt z^3A^3\t)\alpha y^2\t w^2 Y^1R\bt x^2T^3\o X^3\o  \\
	&		& \tens Q^1(y^3\o w^3\o Y^2R\bo x^1 T^2X^2\t R'\bo\la z^1A^2b\t g^2S(z^2A^3\o))S(Q^2)\\
	&		& \,\,\,\,\, y^3\t w^3\t Y^3x^3T^3\t X^3\t\\
	&	=	& X^1A^1b\o g^1S(w^1T^1X^2\o R'\bt z^3A^3\t)\alpha w^2 Y^1R\bt x^2T^3\o X^3\o \\
	&		& \tens Q^1(w^3\o Y^2R\bo x^1 T^2X^2\t R'\bo\la z^1A^2b\t g^2S(z^2A^3\o))S(Q^2)\\
	&		& \,\,\,\,\, w^3\t Y^3x^3T^3\t X^3\t\\
	&	=	& X^1A^1b\o g^1S(w^1T^1\underline{X^2\o R'\bt} z^3A^3\t)\alpha w^2 Y^1R\bt x^2T^3\o X^3\o  \\
	&		& \tens Q^1(w^3\o Y^2R\bo x^1 T^2\underline{X^2\t R'\bo}\la z^1A^2b\t g^2S(z^2A^3\o))S(Q^2)w^3\t Y^3x^3T^3\t X^3\t\\
	&	=	& X^1A^1b\o g^1S(\underline{w^1}T^1R'\bt X^2\t z^3A^3\t)\alpha \underline{w^2 Y^1}R\bt x^2T^3\o X^3\o \\
	&		& \tens Q^1(\underline{w^3\o Y^2}R\bo x^1 T^2R'\bo X^2\o\la z^1A^2b\t g^2S(z^2A^3\o))S(Q^2)\underline{w^3\t Y^3}x^3T^3\t X^3\t\\
	&	=	& X^1A^1b\o g^1S(w^1t^1T^1R'\bt X^2\t z^3A^3\t)\alpha w^2\underline{t^2\o R\bt}x^2T^3\o X^3\o \\
	&		& \tens Q^1(w^3\underline{t^2\t R\bo} x^1 T^2R'\bo X^2\o\la z^1A^2b\t g^2S(z^2A^3\o))S(Q^2)t^3x^3T^3\t X^3\t\\
	&	=	& X^1A^1b\o g^1S(w^1\underline{t^1T^1}R'\bt X^2\t z^3A^3\t)\alpha w^2R\bt \underline{t^2\t x^2T^3\o} X^3\o \\
	&		& \tens Q^1(w^3R\bo \underline{t^2\o x^1 T^2}R'\bo X^2\o\la z^1A^2b\t g^2S(z^2A^3\o))S(Q^2)\underline{t^3x^3T^3\t} X^3\t\\
	&	=	& X^1A^1b\o g^1S(w^1T^1\underline{x^1\o R'\bt} X^2\t z^3A^3\t)\alpha w^2R\bt T^3x^2X^3\o \\
	&		& \tens Q^1(w^3R\bo T^2\underline{x^1\t R'\bo} X^2\o\la z^1A^2b\t g^2S(z^2A^3\o))S(Q^2)x^3X^3\t\\
	&	=	& X^1A^1b\o g^1S(\underline{w^1T^1R'\bt} x^1\t X^2\t z^3A^3\t)\alpha \underline{w^2R\bt T^3}x^2X^3\o \\
	&		& \tens Q^1(\underline{w^3R\bo T^2R'\bo} x^1\o X^2\o\la z^1A^2b\t g^2S(z^2A^3\o))S(Q^2)x^3X^3\t\\
	&	=	& X^1A^1b\o g^1S(R\bt\o Y^2x^1\t X^2\t z^3A^3\t)\alpha R\bt\t Y^3x^2X^3\o \\
	&		& \tens Q^1(R\bo Y^1x^1\o X^2\o\la z^1A^2b\t g^2S(z^2A^3\o))S(Q^2)x^3X^3\t\\
	&	=	& X^1A^1b\o g^1S(Y^2x^1\t X^2\t z^3A^3\t)\alpha Y^3x^2X^3\o \\
	&		& \tens Q^1(Y^1x^1\o X^2\o\la z^1A^2b\t g^2S(z^2A^3\o))S(Q^2)x^3X^3\t\\
	&	=	& X^1A^1b\o g^1S(Y^2\underline{x^1\t X^2\t z^3}A^3\t)\alpha Y^3x^2X^3\o \\
	&		& \tens W^1Y^1\o \underline{x^1\o\o X^2\o\o z^1}A^2b\t g^2S(W^2Y^1\t \underline{x^1\o\t X^2\o\t z^2}A^3\o)\alpha W^3x^3X^3\t\\
	&	=	& \underline{X^1A^1}b\o g^1S(Y^2z^3x^1\t\t \underline{X^2\t\t A^3\t})\alpha Y^3x^2\underline{X^3\o} \\
	&		& \tens W^1Y^1\o z^1x^1\o \underline{X^2\o A^2}b\t g^2S(W^2Y^1\t z^2x^1\t\o \underline{X^2\t\o A^3\o})\alpha\\
	&		& \,\,\,\,\,  W^3x^3\underline{X^3\t}\\
	&	=	& X^1A^1\o b\o g^1S(\underline{Y^2z^3}x^1\t\t y^2\t X^3\o\t A^2\t)\alpha \underline{Y^3}x^2y^3\o X^3\t\o A^3\o  \\
	&		& \tens W^1\underline{Y^1\o z^1}x^1\o y^1X^2A^1\t b\t g^2S(W^2\underline{Y^1\t z^2}x^1\t\o y^2\o X^3\o\o A^2\o)\alpha\\
	&		& \,\,\,\,\, W^3x^3y^3\t X^3\t\t A^3\t\\
	&	=	& X^1A^1\o b\o g^1S(T^2Y^2\t x^1\t\t y^2\t X^3\o\t A^2\t)\alpha T^3Y^3x^2y^3\o X^3\t\o A^3\o \\
	&		& \tens W^1Y^1x^1\o y^1X^2A^1\t b\t g^2S(W^2T^1Y^2\o x^1\t\o y^2\o X^3\o\o A^2\o)\alpha \\
	&		& \,\,\,\,\, W^3x^3y^3\t X^3\t\t A^3\t\\
	&	=	& X^1A^1\o b\o g^1S(T^2\underline{Y^2\t x^1\t\t y^2\t} X^3\o\t A^2\t)\alpha T^3\underline{Y^3x^2y^3\o} X^3\t\o A^3\o \\
	&		& \tens W^1\underline{Y^1x^1\o y^1}X^2A^1\t b\t g^2S(W^2T^1\underline{Y^2\o x^1\t\o y^2\o} X^3\o\o A^2\o)\alpha\\
	&		& \,\,\,\,\, W^3\underline{x^3y^3\t} X^3\t\t A^3\t\\
	&	=	& X^1A^1\o b\o g^1S(T^2y^2\o\t \underline{x^1\t X^3\o\t} A^2\t)\alpha T^3y^2\t \underline{x^2X^3\t\o} A^3\o \\
	&		& \tens W^1y^1X^2A^1\t b\t g^2S(W^2T^1y^2\o\o \underline{x^1\o X^3\o\o} A^2\o)\alpha W^3y^3\underline{x^3X^3\t\t} A^3\t\\
	&	=	& X^1A^1\o b\o g^1S(\underline{T^2y^2\o\t X^3\o\o\t} x^1\t A^2\t)\alpha \underline{T^3y^2\t X^3\o\t} x^2A^3\o \\
	&		& \tens W^1y^1X^2A^1\t b\t g^2S(W^2\underline{T^1y^2\o\o X^3\o\o\o} x^1\o A^2\o)\alpha W^3y^3X^3\t x^3A^3\t\\
	&	=	& X^1A^1\o b\o g^1S(T^2x^1\t A^2\t)\alpha T^3x^2A^3\o \\
	&		& \tens W^1y^1X^2A^1\t b\t g^2S(W^2y^2 X^3\o T^1x^1\o A^2\o)\alpha W^3y^3X^3\t x^3A^3\t\\
	&	=	& X^1A^1\o b\o g^1S(T^2x^1\t A^2\t)\alpha T^3x^2A^3\o \\
	&		& \tens X^2A^1\t b\t g^2S(X^3\o T^1x^1\o A^2\o)\alpha X^3\t x^3A^3\t\\
	&	=	& A^1\o b\o g^1S(T^2x^1\t A^2\t)\alpha T^3x^2A^3\o \tens A^1\t b\t g^2S(T^1x^1\o A^2\o)\alpha x^3A^3\t\\
	&	=	& X^1\o b\o g^1S(X^2\o)S(Y^2x^1\t)\alpha Y^3x^2X^3\o \\
	&		& \tens X^1\t b\t g^2 S(X^2\o)S(Y^1x^1\o)\alpha x^3X^3\t\\
	&	=	& X^1\o b\o S(X^2)\o g^1\gamma^1 X^3\o \tens X^1\t b\t S(X^2)\t g^2 \gamma^2 X^3\t\\
	&	=	& X^1\o b\o S(X^2)\o \alpha\o X^3\o \tens X^1\t b\t S(X^2)\t \alpha\t X^3\t\\
	&	=	& \Delta(X^1bS(X^2)\alpha X^3)\\
	&	=	& \Delta(q^1 b S(q^2))
\end{eqnarray*}

\endproof

\begin{example}
Recall the structure of the twisted quantum double, $D^{\phi}(G)$, for a finite non-abelian group $G$ from \cite{Dijkgraaf90},

\[ (g \otimes \delta_s)(h \otimes \delta_t) = (gh \otimes \delta_s) \delta_{s,gtg^{-1}} \,\,\, \theta_s(g,h) \]
\[ \eta(1) = (e \otimes 1) \]
\[ \Delta(g \otimes \delta_s) = \sum_{ab=s} (g \otimes \delta_a) \otimes (g \otimes \delta_b) \,\,\, \gamma_g(a,b) \]
\[ \varepsilon(g \otimes \delta_s) = \delta_{s,e} \]
\[ S(g \otimes \delta_s) = g^{-1} \otimes \delta_{s^{-1}} \,\,\, \theta_{s^{-1}}^{-1}(g,g^{-1})\gamma_g^{-1}(s,s^{-1}) \]
\[ \alpha = (e \otimes 1) \]
\[ \beta = \sum_g (e \otimes \delta_g) \,\,\, \phi(g^{-1},g,g^{-1}) \]
\[ \phi_D = \sum_{g,h,k} (e \otimes \delta_g) \otimes (e \otimes \delta_h) \otimes (e \otimes \delta_k) \,\,\, \phi(g,h,k) \]
\[ R = \sum_g (e \otimes \delta_g) \otimes (g \otimes 1) \]

where  for all $g,h,t \in G$,

\[ \theta_s(g,h)=\phi(g,g^{-1}sg,h)\phi^{-1}(s,g,h)\phi^{-1}(g,h,h^{-1}g^{-1}sgh) \]
\[ \gamma_g(a,b)=\phi(a,g,g^{-1}bg)\phi^{-1}(a,b,g)\phi^{-1}(gg^{-1}ag,g^{-1}bg) \]

which further satisfy the following identities.

\[ \theta_s(g,h)\theta_s(gh,k)=\theta_s(g,hk)\theta_{g^{-1}sg}(h,k) \]
\[ \gamma_g(a,b)\gamma_g(ab,c)\phi(a,b,c)=\gamma_g(a,bc)\gamma_g(b,c)\phi(g^{-1}ag,g^{-1}bg,g^{-1}cg) \]
\[ \theta_s(g,h)\theta_t(g,h)\gamma_g(s,t)\gamma_h(g^{-1}sg,s^{-1}tg)=\theta_{st}(g,h)\gamma_{gh}(s,t) \]

The adjoint action of $D^{\phi}(G)$ is given by

\begin{eqnarray*}
(g\tens\delta_s)\la(h\tens\delta_t)	& =	& (ghg^{-1}\tens\delta_{gtg^{-1}}) \, \delta_{s,gth^{-1}t^{-1}hg^{-1}} \\
																		&		&	\,\,\, \gamma_g(gtg^{-1},gt^{-1}g^{-1})\gamma_g^{-1}(gh^{-1}t^{-1}hg^{-1},gh^{-1}thg^{-1})\\
																		&		&	\,\,\, \theta_{gtg^{-1}}(g,h)\theta_{gtg^{-1}}(gh,g^{-1})\theta_{b^{-1}}^{-1}(g,g^{-1})
\end{eqnarray*}

We note that $(e\tens\delta_s)\la(h\tens\delta_t) = (h\tens\delta_t) \, \delta_{s,th^{-1}t^{-1}h}$.

We find the structure of $\underline{D^{\phi}(G)}$ to be

\begin{eqnarray*}
\underline{m}((g\tens\delta_s)\tens(h\tens\delta_t)) & = & (gh\tens\delta_s) \, \delta_{s,gtg^{-1}} \theta_s(g,h)\\
	&		&	\,\,\, \phi(s,g^{-1}s^{-1}g,g^{-1}sg)\phi^{-1}(sg^{-1}s^{-1}g,g^{-1}sg,h^{-1}g^{-1}s^{-1}gh) 
\end{eqnarray*}

\[ \underline{\eta}(1) = \sum_{g\in G}(e\tens\delta_g) \, \phi(g^{-1},g,g^{-1}) \]

\begin{eqnarray*}
\underline{\Delta}(g\tens\delta_s)	&	=	&	\sum_{ab=s} (bgb^{-1}\tens\delta_a)\tens(g\tens\delta_b) \, \gamma_g(a,b)\theta_a^{-1}(bgb^{-1},bg^{-1}b^{-1}g)\phi(s,g^{-1}s^{-1}g,g^{-1}sg)\\
	&		&	\,\,\, \phi^{-1}(a,bg^{-1}b^{-1}a^{-1}bgb^{-1},bg^{-1}b^{-1}abgb^{-1})\phi^{-1}(b,g^{-1}b^{-1}g,g^{-1}bg)\\
	&		& \,\,\, \phi(bg^{-1}b^{-1}g,g^{-1}ag,g^{-1}bg)\phi^{-1}(abg^{-1}b^{-1}a^{-1}bgb^{-1},bg^{-1}b^{-1}g,g^{-1}abg)\\
	&		& \,\,\, \phi(abg^{-1}b^{-1}a^{-1}bgb^{-1},bg^{-1}b^{-1}abgb^{-1},b)\phi^{-1}(bg^{-1}b^{-1}abgb^{-1},bg^{-1}b^{-1}g,g^{-1}bgb)
\end{eqnarray*}

\[ \underline{\varepsilon}(g\tens\delta_s) = \delta_{s,e} \]

\begin{eqnarray*}
\underline{S}(g\tens\delta_s)	&	=	& (sg^{-1}s^{-1}\tens\delta_{sg^{-1}s^{-1}gs^{-1}})\theta_{s^{-1}}^{-1}(g,g^{-1}\gamma_g^{-1}(s,s^{-1})\\
	&		&	\,\,\, \theta_{sg^{-1}s^{-1}gs^{-1}}(sg^{-1}s^{-1}g,g^{-1})\phi(sg^{-1}s^{-1}gs^{-1},sg^{-1}s^{-1}g,g^{-1}sg)\\
	&		& \,\,\, \phi(s,g^{-1}s^{-1}g,g^{-1}sg)\phi^{-1}(sg^{-1}s^{-1}g,g^{-1}s^{-1}g,g^{-1}sg)\phi(g^{-1}sg,g^{-1}s^{-1}g,g^{-1}sg)
\end{eqnarray*}
\end{example}

\section{Bosonisation of Braided Groups in $_H\Mcal$}

Let $H$ be a quasi-triangular quasi-Hopf algebra. Given a braided group in $_{H}\mathcal{M}=\mathcal{C}$ we can `\,bosonise' it back to an equivalent ordinary quasi-Hopf algebra. We use the same strategy as in \cite{Majid94}. If $B$ is a braided Hopf algebra in $\mathcal{C}$, then a braided $B$-module is an object $V$ in $\mathcal{C}$ and a morphism $\alpha^{B}_{V}: B \tens V \rightarrow V$ in $\mathcal{C}$. Note that $\alpha^{B}_{V}$ intertwines the action of $H$, that is $\alpha^{B}_{V}(h \la (b \tens v)) = h \la \alpha^{B}_{V} (b \tens v)$, for all $h \in H, b \in B, v \in V$; equivalently,

\[ h \la (b \la v) = (h\o\la b) \la (h\t \la v) \]

where the notation for the actions of $H$ on $B$, $H$ on $V$ and $B$ on $V$ are understood. The category $_{B}\mathcal{C}$ of braided $B$-modules in $\mathcal{C}$ is a braided monoidal category with the same braiding as $\mathcal{C}$.

\begin{theorem} Let $H$ be a quasitriangular quasi-Hopf algebra, and $B \in \mathcal{C}$ be a braided group. Then there is an ordinary quasi-Hopf algebra $B \bos H$ built on the vector space $B \tens H$ with structure

\[ (b \tens h)(c \tens g) = (x^1 \la b)\underline{\cdot}(x^2h\o \la c) \tens x^3h\t g \]
\[ \eta(1) = 1_B \tens 1 \]
\[ \Delta(b \tens h) = y^{1}X^{1} \la b\underline{\o} \tens y^{2}Y^{1}R\bt x^{2}X^{3}\o h\o \tens y^{3}\o Y^{2}R\bo x^{1}X^{2} \la b\underline{\t} \tens y^{3}\t Y^{3}x^{3}X^{3}\t h\t \]
\[ \varepsilon(b \tens h) = \underline{\varepsilon}(b)\varepsilon(h) \]
\[ S(b \tens h) =  (S(X^1x^1\o R\bt h)\alpha)\o X^2x^1\t R\bo \la \underline{S}(b) \tens (S(X^1x^1\o R\bt h)\alpha)\t X^3x^2\beta S(x^3) \]
\[ \alpha_{B \rtimes\kern -3.75pt \cdot H} = 1_B \tens \alpha \]
\[ \beta_{B \rtimes\kern -3.75pt \cdot H} = 1_B \tens \beta \]
\[ \phi_{B \rtimes\kern -3.75pt \cdot H} = 1_B \tens X^1 \tens 1_B \tens X^2 \tens 1_B \tens X^3 \]
\end{theorem}

\proof
Given a braided $B$-module, $V$, in $\mathcal{C}$, we have an action of $B$ on $V$ and an action of $H$ on $V$. The action of $B \bos H$ on $V$ is

\[ (b \tens h) \la v = b \la (h \la v) \]

for all $v \in V, b \in B, h \in H$. Note, that since the action of $B$ on $V$ is a morphism in $\mathcal{C}$ it satisfies

\[ b \la (c \la v) = (x^1 \la b)(x^2 \la c) \la (x^3 \la v) \]

for all $b,c \in B, v \in V$. Since $B \bos H$ is an ordinary Hopf algebra, it satisfies

\[ (b \tens h)(c \tens g) \la v = (b \tens h) \la ((c \tens g) \la v) \]

and hence this determines the multiplication in $B \bos H$.

\begin{eqnarray*}
\lefteqn{ (b \tens h)(c \tens g) \la v = b \la (h \la ((c \tens g) \la v)) } \\
 & = & b \la (h \la (c \la (g \la v))) \\
 & = & b \la ((h\o \la c) \la (h\t \la (g \la v))) \\
 & = & b \la ((h\o \la c) \la (h\t g \la v)) \\
 & = & (x^1 \la b)(x^2 \la (h\o \la c)) \la (x^3 \la (h\t g \la v)) \\
 & = & (x^1 \la b)(x^2h\o \la c) \la (x^3h\t g \la v) 
\end{eqnarray*}

Thus,

\[ (b \tens h)(c \tens g) = (x^1 \la b)(x^2h\o \la c) \tens x^3h\t g. \]

Let $V,W \in {}_B\mathcal{C}$then $B \bos H$ acts on $V \tens W$ as

\[ (b \tens h) \la (v \tens w) = (b \tens h)\o \la v \tens (b \tens h)\t \la w \]

for all $v \in V, w \in W$. But also

\[ (b\tens h) \la (v \tens w) = b \la (h \la (v \tens w)) = b \la (h\o \la v \tens h\t \la w) \]

Thus the coproduct of $B \bos H$ is characterised by

\[ (b \tens h)\o \la v \tens (b \tens h)\t \la w = b \la (h\o \la v \tens h\t \la w) \]

Now, $B$ acts on the tensor product $V \tens W$ as

\[ \alpha^B_{V \tens W} = (\alpha^B_V \tens \alpha^B_W)  \Phi^{-1}_{B,V,B \tens W}  (\id \tens \Phi_{V,B,W})  (\id \tens \Psi_{B,V} \tens \id)  (\id \tens \Phi^{-1}_{B,V,W})  \Phi_{B,B,V \tens W}  (\underline{\Delta} \tens \id \tens \id) \]

that is,

\begin{eqnarray*}
b \la (v \tens w) &	=	& (y^1X^1 \la b\underline{\o}) \la (y^2Y^1R\bt x^2X^3\o \la v) \tens (y^3\o Y^2R\bo x^1X^2 \la b\underline{\t}) \la (y^3\t Y^3x^3X^3\t \la w)
\end{eqnarray*}

So,

\begin{eqnarray*}
\lefteqn{(b \tens h) \la (v \tens w) = b \la (h\o \la v \tens h\t \la w)} \\
 & = & (y^1X^1 \la b\underline{\o}) \la (y^2Y^1R\bt x^2X^3\o h\o \la v) \tens (y^3\o Y^2R\bo x^1X^2 \la b\underline{\t}) \la (y^3\t Y^3x^3X^3\t h\t \la w)
\end{eqnarray*}

Thus,

\[ \Delta(b \tens h) = y^1X^1 \la b\underline{\o} \tens y^2Y^1R\bt x^2X^3\o h\o \tens y^3\o Y^2R\bo x^1X^2 \la b\underline{\t} \tens y^3\t Y^3x^3X^3\t h\t \]

For the antipode, given $V \in {}_B\mathcal{C}$ we have to consider how $H$ and $B$ act on the dual object $V^*$. It is known that for a quasi-Hopf algebra $H$ and a left $H$-module $V$, the dual space $V^*$ becomes a left $H$-module by $(h \la v^*)(v) = v^*(h \la v)$ for any $v^* \in V^*, v \in  V, h \in H$, thus for $B \bos H$,

\[((b \tens h) \la v^*)(v) = v^*(S(b \tens h) \la v) \]

for all $v^* \in V^*, v \in V, b \in B, h \in H$. But we also have

\[ ((b \tens h) \la v^*)(v) = (b \la (h \la v^*))(v) \]

so the antipode is determined by

\[ v^*(S(b \tens h) \la v) = (b \la (h \la v^*))(v) \]

so it remains to find how $B$ acts on the dual space. If $V$ is a left $B$-module, then $V^*$ is a right $B$-module by $\alpha^*: V^* \tens B \rightarrow V^*$ as

\[ \alpha^* = l^{-1}_{V^*}  (ev_V \tens \id)  (\id \tens \alpha^B_V \tens \id)  (\id \tens \id \tens coev_V) (\id \tens r_B) \]

so,

\begin{eqnarray*}
\lefteqn{(\alpha^*(v^* \tens b))(v) = (l^{-1}_{V^*}  \ldots  (\id \tens r_B)(v^* \tens b))(v)} \\
 & = & l^{-1}_{V^*}  \ldots  \Phi^{-1}_{V^*,B,\underline{1}} (v^* \tens (b \tens 1))(v) \\
 & = & l^{-1}_{V^*}  \ldots  (\id \tens \id \tens coev_V)((x^1 \la v^* \tens x^2 \la b) \tens x^3 \la 1)(v) \\
 & = & l^{-1}_{V^*}  \ldots  (\id \tens \id \tens coev_V)((x^1 \la v^* \tens x^2 \la b) \tens \varepsilon(x^3))(v) \\
 & = & l^{-1}_{V^*}  \ldots  (\id \tens \id \tens coev_V)((v^* \tens b) \tens 1)(v) \\
 & = & l^{-1}_{V^*}  \ldots  \Phi^{-1}_{V^* \tens B,V,V^*}((v^* \tens b) \tens (\beta \la e_a \tens f^a))(v) \\
 & = & l^{-1}_{V^*}  \ldots  (\Phi_{V^*,B,V} \tens \id)(((x^1\o \la v^* \tens x^1\t \la b) \tens x^3\beta \la e_a) \tens x^3 \la f^a)(v) \\
 & = & l^{-1}_{V^*}  \ldots  (\id \tens \alpha^B_V \tens \id)((X^1x^1\o \la v^* \tens (X^2x^1\t \la b \tens X^3x^2\beta \la e_a)) \tens x^3 \la f^a)(v) \\
 & = & l^{-1}_{V^*}  (ev_V \tens \id)((X^1x^1\o \la v^* \tens (X^2x^1\t \la b) \la (X^3x^2\beta \la e_a)) \tens x^3 \la f^a)(v) \\
 & = & l^{-1}_{V^*}((X^1x^1\o \la v^*)(\alpha \la ((X^2x^1\t \la b) \la (X^3x^2\beta \la e_a))) \tens x^3 \la f^a)(v) \\
 & = & l^{-1}_{V^*}(v^*(S(X^1x^1\o)\alpha \la ((X^2x^1\t \la b) \la (X^3x^2\beta \la e_a))) \tens x^3 \la f^a)(v) \\
 & = & v^*(S(X^1x^1\o)\alpha \la ((X^2x^1\t \la b) \la (X^3x^2\beta S(x^3) \la v)))
\end{eqnarray*}

Then, $V^*$ becomes a left $B$-module by

\begin{eqnarray*}
\lefteqn{\alpha^B_{V^*}(b \tens v^*)(v) = \alpha^*  (\id \tens \underline{S})  \Psi_{B,V^*}(b \tens v^*)(v)} \\
 & = & \alpha^*  (\id \tens \underline{S})(R\bt \la v^* \tens R\bo \la b)(v) \\
 & = & \alpha^* (R\bt \la v^* \tens \underline{S}(R\bo \la b))(v) \\
 & = & \alpha^* (R\bt \la v^* \tens R\bo \la \underline{S}(b))(v) \\
 & = & (R\bt \la v^*)(S(X^1x^1\o)\alpha \la ((X^2x^1\t R\bo \la \underline{S}(b)) \la (X^3x^2\beta S(x^3) \la v))) \\
 & = & v^*(S(R\bt)S(X^1x^1\o)\alpha \la ((X^2x^1\t R\bo \la \underline{S}(b)) \la (X^3x^2\beta S(x^3) \la v)))
\end{eqnarray*}

So, the action of $B \bos H$ on $V^*$ is given by

\begin{eqnarray*}
\lefteqn{((b \tens h) \la v^*)(v) = (b \la (h \la v^*))(v)} \\
 & = & (h \la v^*)(S(R\bt)S(X^1x^1\o)\alpha \la ((X^2x^1\t R\bo \la \underline{S}(b)) \la (X^3x^2\beta S(x^3) \la v))) \\
 & = & v^*(S(h)S(R\bt)S(X^1x^1\o)\alpha \la ((X^2x^1\t R\bo \la \underline{S}(b)) \la (X^3x^2\beta S(x^3) \la v))) \\
 & = & v^*(S(X^1x^1\o R\bt h)\alpha \la ((X^2x^1\t R\bo \la \underline{S}(b)) \la (X^3x^2\beta S(x^3) \la v))) \\
 & = & v^*(((S(X^1x^1\o R\bt h)\alpha)\o X^2x^1\t R\bo \la \underline{S}(b)) \la ((S(X^1x^1\o R\bt h)\alpha)\t X^3x^2\beta S(x^3) \la v))
\end{eqnarray*}

So

\begin{eqnarray*}
\lefteqn{((b \tens h) \la v^*)(v) = v^*(S(b \tens h) \la v)} \\
 & = & v^*(((S(X^1x^1\o R\bt h)\alpha)\o X^2x^1\t R\bo \la \underline{S}(b)) \la ((S(X^1x^1\o R\bt h)\alpha)\t X^3x^2\beta S(x^3) \la v))
\end{eqnarray*}

Hence,

\[ S(b \tens h) = (S(X^1x^1\o R\bt h)\alpha)\o X^2x^1\t R\bo \la \underline{S}(b) \tens (S(X^1x^1\o R\bt h)\alpha)\t X^3x^2\beta S(x^3). \]

\endproof

\begin{corollary} The modules of $B$ in the braided category $_H\mathcal{M}$ correspond to the ordinary modules of $B \bos H$. The braided categories ar isomorphic.
\end{corollary}

\proof $B\bos H$ is a smash product when considered as an algebra. This structure was found in \cite{Bulacu00}, and this correspondence is given as follows. Let $V$ be a $B \bos H$-module with structure given by $(b \tens h) \cdot v$. Define maps $j: H \rightarrow B \bos H$ and $i: B \rightarrow B \bos H$ by $j(h) = 1 \tens h$ and $i(b) = b \tens 1$. Then $V$ becomes a left $H$-module by $h \la v = j(h) \cdot v$, and $V$ becomes a braided $B$-module by $b \la v = i(b) \cdot v$. Conversely, if $V$ is a braided module in $_H\mathcal{M}$, define the action of $B \bos H$ on $V$ by $(b \tens h) \cdot v = b \la (h \la v)$. It is straightforward to see that this is an equivalence of monoidal categories by the same steps as in \cite{Majid94}.

\endproof

\begin{example}For an example of braided group bosonisation, we consider the group function algebra $k_{\phi}(G)$, and find an isomorphism $\underline{kG}\bos k_{\phi}(G) \cong D^{\phi}(G)$. 

Consider the group function algebra, $k(G)$, of a finite group $G$ with identity $e$. This is the set of functions on $G$ with values in $k$. This has the structure of a commutative Hopf algebra as follows.

\[ \delta_s \cdot \delta_t = \delta_t \,\,\, \delta_{s,t} \]
\[ \eta(1) = \sum_{t \in G} \delta_t \]
\[ \Delta(\delta_t) = \sum_{ab=t} \delta_a \tens \delta_b \]
\[ \varepsilon(\delta_t) = \delta_{t,e} \]
\[ S(\delta_t) = \delta_{t^{-1}} \]

for all $\delta_s,\delta_t \in k(G)$. Further, one can view $k(G)$ as a quasi-Hopf algebra with $\phi_G \in k(G)^{\tens 3}$ defined by

\[ \phi_G = \sum_{r,s,t \in G} \delta_r \tens \delta_s \tens \delta_t \,\,\, \phi(r,s,t) \]

for some 3-cocycle $\phi \in k(G)$ satisfying

\[ \phi(b,c,d)\phi(a,bc,d)\phi(a,b,c) = \phi(a,b,cd)\phi(ab,c,d) \]
\[ \phi(a,e,b) = 1 \]

for all $a,b,c,d \in G$. Choosing $\alpha = \varepsilon_{kG} = 1$, this determines $\beta \in k(G)$ as

\[ \beta =  \sum_{t \in G}\delta_t \,\,\, \phi^{-1}(t,t^{-1},t) = \sum_{t \in G}\delta_t \,\,\, \phi(t^{-1},t,t^{-1}) \]

A quasitriangular structure for $k(G)$ as a quasi-Hopf algebra is defined by $R= \sum_{s,t \in G} \delta_s \tens \delta_t \,\,\, r(s,t)$, where $r \in k(G) \tens k(G)$  is a function obeying

\[ r(gh,t) = r(g,t)r(h,t)\frac{\phi(t,g,h)\phi(g,h,t)}{\phi(g,t,h)} \]
\[ r(t,gh) = r(t,g)r(t,h)\frac{\phi(g,t,h)}{\phi(t,g,h)\phi(g,h,t)} \]
\[ r(u,e) = 1 = r(e,u) \]

for all $g,h,t \in G$. We denote this quasitriangular quasi-Hopf algebra by $k_{\phi}(G)$. The structure of $\underline{k_{\phi}(G)}$ is as follows

\[ \underline{m}(\delta_s\tens\delta_t) = \delta_t \, \delta_{s,t} \phi(t,t^{-1},t) \]

\[ \underline{\eta}(1) =  \sum_{s \in G}\delta_s \,\,\, \phi(s^{-1},s,s^{-1})\]

\[ \underline{\Delta}(\delta_t) = \sum_{ab=t} \delta_a \tens \delta_b \,\,\, \frac{\phi(t,t^{-1},t)}{\phi(a,a^{-1},a)\phi(b,b^{-1},b)} \]

\[ \underline{\varepsilon}(\delta_s) = \delta_{s,e} \]

\[ \underline{S}(\delta_s) = \delta_{s^{-1}} \,\,\, \phi(s,s^{-1},s)\phi(s,s^{-1},s) \]

So $\underline{k(G)}$ has structure

\[ \delta_s \underline{\cdot} \delta_t = \delta_t \,\,\, \delta_{s,t} \phi(t,t^{-1},t) \]
\[ \underline{\eta}(1)= \sum_{t \in G} \delta_t \phi(t^{-1},t,t^{-1}) \]
\[ \underline{\Delta}\delta_t = \sum_{ab=t} \delta_a \tens \delta_b \,\,\, \frac{\phi(t,t^{-1},t)}{\phi(a,a^{-1},a)\phi(b,b^{-1},b)}  \]
\[ \underline{\varepsilon}(\delta_t) = \delta_{t,e} \]
\[ \underline{S}(\delta_t) = \delta_{t^{-1}} \,\,\, \phi(t,t^{-1},t)\phi(t,t^{-1},t) \]

We can find the structure on $\underline{kG} = \underline{(k_{\phi}(G))^*} = (\underline{k_{\phi}(G)})^*$ as this braided dual structure is determined by the structure on the original braided group. First, consider how $\underline{kG} \in {}_{k_{\phi}(G)}\mathcal{M}$. Let $\psi \in k_{\phi}(G), \delta_s \in \underline{k_{\phi}(G)}, g \in \underline{kG}$, then

\[ \langle \psi \la g, \delta_s \rangle = \langle g, S(\psi) \la \delta_s \rangle \]

So,

\begin{eqnarray*}
\delta_s(\psi \la g) & = & (S(\psi) \la \delta_s)(g) \\
 & = & ((S\psi)\o\delta_s S((S\psi)\t))(g) \\
 & = & (S\psi)\o(g) \delta_s(g) (S\psi)\t(g^{-1}) \\
 & = & (S\psi)(e)\delta_s(g) \\
 & = & \psi(e) \delta_s(g) \\
 & = & \delta_s(\psi(e)g)
\end{eqnarray*}

Thus, $\psi \la g = \psi(e)g$ for all $\psi \in k_{\phi}(G), g \in \underline{kG}$. So, if we consider the associativity constraint $\Phi$ on this category; if it is acting on $\underline{kG}$, it is in fact trivial, and as such, in the following calculations we can ignore the bracketing order. Now, the multiplication on $\underline{kG}$ is determined by the comultiplication on $\underline{k_{\phi}(G)}$ as follows:

\[ ev  (r^{-1} \tens \id)  (\id \tens ev \tens \id)  (\underline{\Delta} \tens \id \tens \id)= ev  (\id \tens \underline{m}) : \underline{k(G)} \tens \underline{kG} \tens \underline{kG} \rightarrow \underline{1} \]

The left hand side gives

\begin{eqnarray*}
\lefteqn{ev  (r^{-1} \tens \id)  (\id \tens ev \tens \id)  (\underline{\Delta} \tens \id \tens \id)(\delta_s \tens g \tens h)} \\
 & = & ev  (r^{-1} \tens \id)  (\id \tens ev \tens \id)((\delta_s)\underline{\o} \tens (\delta_s)\underline{\t} \tens g \tens h) \\
 & = & ev ((\delta_s)\underline{\o} \tens h \, (\delta_s)\underline{\t}(\alpha \la g)) \\
 & = & (\delta_s)\underline{\o}(\alpha \la h)(\delta_s)\underline{\t}(\alpha \la g) \\
 & = & \underline{\Delta}(\delta_s)(h,g) \\
 & = & \frac{\phi(hg,(hg)^{-1},hg)}{\phi(h,h^{-1},h)\phi(g,g^{-1},g)}\delta_s(hg) \\
 & = & \frac{\phi(gh,(gh)^{-1},gh)}{\phi(g,g^{-1},g)\phi(h,h^{-1},h)}\delta_s(gh)
\end{eqnarray*}

while the right hand side gives

\begin{eqnarray*}
ev  (\id \tens \underline{m})(\delta_s \tens g \tens h) & = & ev(\delta_s \tens g \underline{\cdot} h) \\
 & = & \delta_s(\alpha \la (g \underline{\cdot} h)) \\
 & = & \delta_s(g \underline{\cdot} h)
\end{eqnarray*}

These are equal, hence $\underline{m}(g \tens h) = \frac{\phi(gh,(gh)^{-1},gh)}{\phi(g,g^{-1},g)\phi(h,h^{-1},h)} \, gh$ for all $g,h \in \underline{kG}$. The rest of the structure is similarly determined; the structure of $\underline{kG}$ is

\[ \underline{m}(g \tens h) = gh \,\,\, \frac{\phi(gh,(gh)^{-1},gh)}{\phi(g,g^{-1},g)\phi(h,h^{-1},h)} \]
\[ \underline{\eta}(1) = e \]
\[ \underline{\Delta}(g) = g \tens g \,\,\, \phi(g,g^{-1},g) \]
\[ \underline{\varepsilon}(g) = \phi(g^{-1},g,g^{-1}) \]
\[ \underline{S}(g) = g^{-1} \,\,\, \phi(g^{-1},g,g^{-1})\phi(g^{-1},g,g^{-1}) \]

Finally, we can bosonise $\underline{kG} \in {}_{k_{\phi}(G)}\mathcal{M}$ into an ordinary quasi-Hopf algebra with the following structure:

\[ (g \tens \delta_s)(h \tens \delta_t) = gh \tens \delta_{s,t}\delta_t \frac{\phi(gh,(gh)^{-1},gh)}{\phi(g,g^{-1},g)\phi(h,h^{-1},h)} \]
\[ \eta(1)(g \tens \delta_t) = e \tens 1 \]
\[ \Delta(g \tens \delta_t) = \sum_{ab = t} g \tens \delta_a \tens g \tens \delta_b \, \phi(g,g^{-1},g) \]
\[ \varepsilon(g \tens \delta_t) = \phi(g^{-1},g,g^{-1})\delta_{t,e} \]
\[ S(g \tens \delta_t) = g^{-1} \tens \delta_{t^{-1}} \phi(g^{-1},g,g^{-1})\phi(g^{-1},g,g^{-1}) \]

There exists a quasi-Hopf algebra isomorphism $\sigma: \underline{kG} \bos k_{\phi}(G) \rightarrow D^{\phi}(G)$defined by

\[ \sigma (g \tens \delta_t) = g \tens \delta_t \,\,\, \phi(g^{-1},g,g^{-1})R^{-1}(g,t) \]

It is straightforward to check that $\sigma$ is an isomorphism of quasi-Hopf algebras. Using this isomorphism and its inverse, one can obtain the quasitriangular structure of $\underline{kG} \bos k_{\phi}(G)$. Note that $\sigma^{-1}(g \tens \delta_t) = g \tens \delta_t \,\,\, \phi^{-1}(g^{-1},g,g^{-1})R(g,t) = g \tens \delta_t \,\,\, \phi(g,g^{-1},g)R(g,t)$, hence,

\begin{eqnarray*}
R_B & = & (\sigma^{-1} \tens \sigma^{-1})(R_D) \\
 & = & \sum_{g \in G}(\sigma^{-1} \tens \sigma^{-1})(e \tens \delta_g \tens g \tens 1) \\
 & = & \sum_{g,h \in G} e \tens \delta_g \tens g \tens 1 \,\,\, \phi(g,g^{-1},g)R(g,h) \\
 & = & \sum_{g \in G} e \tens \delta_g \tens g \tens 1 \,\,\, \phi(g,g^{-1},g)R(g,g)
\end{eqnarray*}

\end{example}

\begin{remark}
Similarly, we would expect $\underline{H}^*\bos H \cong D^{\phi}(H)$ for any quasi-Hopf algebra $H$. When $H$ is factorisable, $\underline{H}^*\cong \underline{H}$, and this case is covered in the next section.
\end{remark}

\begin{example} Following \cite{Majid00}, we consider a group $G$ and an invetible 2-cochain $F: G \times G \rightarrow k^*$ satisfying $F(e,g)=F(g,e)=1$ for all $g \in G$. Then one can consider the deformation of the group algebra $kG$ with modified product

\[ g \cdot_F h = F(g,h)gh \]

for all $g,h \in G$ and where $gh$ is the usual group product in $G$.

For a group $G$ and $\phi:G \times G \times G \rightarrow k^*$ an invertible group 3-cocycle, the category of $G$-graded vector spaces (the category of $k(G)$-modules is monoidal with associator determined by $\phi$ and the grading. From \cite{Majid99}, $k_FG$ is a $G$-graded quasialgebra with $|g|=g$ for $g \in G$, which is quasiassociative with associator $\phi$ the coboundary of $F$, that is,

\[ (g \cdot_F h) \cdot_F k = \phi(|g|,|h|,|k|) g \cdot_F (h \cdot_F k) \]
\[ \phi(g,h,k) = \frac{F(g,h)F(gh,k)}{F(h,k)F(g,hk)} \]

for all $g,h,k \in G$. Here $\phi: G \times G \times G \rightarrow k^*$ is an invertible 3-cocycle, and gives the category of $G$-graded vector spaces a monoidal structure.

If $G$ is abelian and $\phi$ is of coboundary form ($\phi = \partial F$), the category of $G$-graded spaces is braided with $\Psi$ determined by the function $R(g,h) = \frac{F(g,h)}{F(h,g)}$ and $kG_F$ is quasicommutative with $g \cdot_F h = R(g,h) h \cdot_F g$.

In the case when $G = \mathbb{Z}_2^n$, $F$ takes the form $F(g,h)=(-1)^{f(g,h)}$ where $f$ is a $\mathbb{Z}_2$-valued function on $G \times G$ such that $F^2=1$. The octonions are of this form for the group group $G = \mathbb{Z}_2 \times \mathbb{Z}_2 \times \mathbb{Z}_2$ with

\[ f(g,h) = \sum_{i \leq j} g_ih_j + h_1g_2g_3 + g_1h_2g_3 + g_1g_2h_3 \]
\[ \phi(g,h,k) = (-1)^{(g \times h)\cdot k} = (-1)^{|ghk|} \]
\begin{equation*}
R(g,h) = \left\{
	\begin{array}{rl}
		1 	& \text{if } g=e \text{ or } h=e \text{ or } g=h \\
		-1 	& \text{otherwise}
	\end{array} \right.
\end{equation*}

We let $G=\mathbb{Z}_2^3$, and consider the graded basis $\{e_a|a\in \mathbb{Z}_2\}$ of $kG_F$ and the dual basis of delta functions $\{\delta_a \}$ of the group function algebra $k(G)$. We can view a $G$-graded quasialgebra as a $k_{\phi}(G)$-module quasi-algebra with action $\delta_b \la e_a = \delta_b(|e_a|)e_a=\delta_{b,a}e_a$ on homogeneous elements, where $k_{\phi}(G)$ is the usual group function algebra on $G$ regarded as a quasi-Hopf algebra with $\phi \in k(G)^{\tens 3}$. Thus, the algebra of octonions,$kG_F$, live naturally in the category of $k_{\phi}(G)$-modules, and as such we can construct the bosonisation of the octonions as an algebra.

\begin{eqnarray*}
(e_a \tens \delta_s)(e_b \tens \delta_t) & = & (\phi^{-(1)} \la e_a)\cdot_F(\phi^{-(2)}(\delta_s)\o \la e_b) \tens \phi^{-(3)}(\delta_s)\t\delta_t \\
 & = & \sum_{xy=s} (\phi^{-(1)} \la e_a)\cdot_F(\phi^{-(2)}\delta_x \la e_b) \tens \phi^{-(3)}\delta_y\delta_t \\
 & = & \sum_{xy=s} \phi^{-(1)}(|e_a|)\phi^{-(2)}(|e_b|)\delta_a(|e_b|) \, e_a \cdot_F e_b \tens \phi^{-(3)}\delta_y\delta_t \\
 & = & \sum_{xy=s} e_a\cdot_F e_b\tens \phi^{-(1)}(a)\phi^{-(2)}(b)\delta_x(b)\phi^{-(3)}\delta_b\delta_t \\
 & = & e_a \cdot_F e_b \tens \delta_t \delta_{-b+s,t}\phi^{-1}(a,b,t) \\
 & = & (-1)^{|abt|} \, e_a\cdot_F e_b \tens \delta_{-b+s,t}\delta_t
\end{eqnarray*}

for all $a,b,s,t \in \mathbb{Z}_2^3$.

It is clear that $(1\tens \delta_s)(1\tens \delta_t)=1\tens \delta_s\delta_t$, and so $\mathbb{O}\bos k_{\phi}(\mathbb{Z}_2^3)\supset k_{\phi}(\mathbb{Z}_2^3)$ as a subalgebra. It also contains an algebra with the following structure.

\begin{eqnarray*}
(e_a\tens 1)(e_b\tens 1)	&	=	& \sum_{s,t} (-1)^{|abt|}e_a \cdot_F e_b \tens \delta_{-b+s,t}\\
	&	=	&	\sum_t (-1)^{|abt|} e_a \cdot_F e_b \tens \delta_t\\
	&	=	& e_a\cdot_F e_b \chi(a,b)
\end{eqnarray*}

where $\chi(a,b)=\sum_t(-1)^{|abt|}\delta_t$, and we note

\begin{equation*}
\chi(a,b) = \left\{
	\begin{array}{ll}
		1 & \text{if } a=0 \text{ or } b=0 \text{ or } a=b \\
		2(\delta_0+\delta_a+\delta_b+\delta_{a+b})-1 	& \text{otherwise}
	\end{array} \right.
\end{equation*}

Finally, the commutation relations are

\begin{eqnarray*}
(e_a\tens 1)(1\tens \delta_t)	&	=	& \sum_s e_a\tens \delta_{s,t}\delta_t\\
	&	=	& e_a\tens \delta_t
\end{eqnarray*}

\begin{eqnarray*}
(1\tens \delta_s)(e_b\tens 1)	&	=	& \sum_t e_b\tens \delta_{s-b,t}\delta_t\\
	&	=	& e_b\tens \delta_{s-b}
\end{eqnarray*}

So we find that $fe_a = e_aL_a(f)$ for all $f\in k_{\phi}(\mathbb{Z}_2^3)$ and $a\in \mathbb{Z}_2^3$, where $L_a(f)(s) = f(a+s)$.

\end{example}

\section{An Isomorphism $\underline{H} \bos H \cong H _{\Rcal}\crosscop H$}

Let $H _{\Rcal}\crosscop H$ be the quasi-Hopf algebra with tensor product algebra, and coproduct
\begin{eqnarray*}
\Delta (b\tens h) 	&	=	&	x^1Y^1 b\o y^1X^1 \tens x^2T^1R\bt w^2Y^3\o h\o y^3\o W^2R^-\bt t^1X^2  \\
	&		& \,\,\,\,\, \tens x^3\o T^2R\bo w^1Y^2 b\t y^2W^1R^-\bo t^2X^3\o\tens x^3\t T^3w^3Y^3\t h\t y^3\t W^3t^3X^3\t
\end{eqnarray*}

\begin{theorem} Let $H$ be a quasitriangular quasi-Hopf algebra. There is a quasi-Hopf algebra isomorphism $\underline{H}\bos H \cong H_\Rcal\crosscop H$ defined by
\[ \chi(a \tens h) = q^1(x^1 \la a)S(q^2)x^2h\o \tens x^3h\t \]
\end{theorem}

\proof
It is straightforward to check that the inverse map is $\chi^{-1}(a \tens h) = x^1aX^1\beta S(x^2h\o X^2) \tens x^3h\t X^3$. First we show that $\chi$ is an algebra morphism,

\begin{eqnarray*}
\chi((a\tens h)(b\tens g))	&	=	&	\chi(q^1(y^1x^1\la a)S(q^2)y^2(x^2h\o\la b)S(y^3) \tens x^3h\t g)\\
	&	=	&	Q^1(w^1\la(q^1(y^1x^1\la a)S(q^2)y^2(x^2h\o\la b)S(y^3)))S(Q^2)w^2x^3\o h\t\o g\o \\
	&		&	\,\,\,\, \tens w^3x^3\o h\t\t g\t\\
	&	=	&	\underline{X^1w^1\o}Y^1y^1\o x^1\o aS(Y^2y^1\t x^1\t)\alpha Y^3y^2x^2\o h\o\o b\\
	&		& \,\,\,\, S(\underline{X^2w^1\t} y^3 x^2\t h\o\t)\alpha \underline{X^3w^2}x^3\o h\t\o g\o \\
	&		&	\,\,\,\, \tens \underline{w^3}x^3\t h\t\t g\t\\
	&	=	&	X^1\underline{Y^1y^1\o} x^1\o aS(\underline{Y^2y^1\t} x^1\t)\alpha \underline{Y^3y^2}x^2\o h\o\o b\\
	&		& \,\,\,\, S(w^1X^2 \underline{y^3} x^2\t h\o\t)\alpha w^2X^3\o x^3\o h\t\o g\o \\
	&		&	\,\,\,\, \tens w^3X^3\t x^3\t h\t\t g\t\\
	&	=	& X^1\underline{Y^1x^1\o} a S(y^1\underline{Y^2x^1\t})\alpha y^2\underline{Y^3\o x^2\o} h\o\o b\\
	&		& \,\,\,\, S(w^1X^2y^3\underline{Y^3\t x^2\t} h\o\t)\alpha w^2X^3\o \underline{x^3\o} h\t\o g\o\\
	&		&	\,\,\,\, \tens w^3X^3\t \underline{x^3\t} h\t\t g\t\\
	&	=	& X^1t^1Y^1aS(\underline{y^1t^2\o} x^1Y^2)\alpha \underline{y^2t^2\t\o} x^2\o Y^3\o\o h\o\o b\\
	&		& \,\,\,\, S(w^1X^2\underline{y^3t^2\t\t} x^2\t Y^3\o\t h\o\t) \alpha w^2X^3\o t^3\o x^3\o \\
	&		& \,\,\,\, Y^3\t\o h\t\o g\o\\
	&		&	\,\,\,\, \tens w^3X^3\t t^3\t x^3\t Y^3\t\t h\t\t g\t\\
	&	=	& \underline{X^1t^1}Y^1aS(y^1x^1Y^2)\alpha y^2x^2\o Y^3\o\o h\o\o b\\
	&		& \,\,\,\, S(w^1\underline{X^2t^2}y^3x^2\t Y^3\o\t h\o\t)\alpha w^2\underline{X^3\o t^3\o} x^3\o \\
	&		& \,\,\,\, Y^3\t\o h\t\o g\o\\
	&		& \,\,\,\, \tens w^3\underline{X^3\t t^3\t} x^3\t Y^3\t\t h\t\t g\t\\
	&	=	& Y^1aS(\underline{y^1x^1}Y^2)\alpha \underline{y^2x^2\o} Y^3\o\o h\o\o b\\
	&		& \,\,\,\, S(w^1\underline{y^3x^2\t} Y^3\o\t h\o\t)\alpha w^2\underline{x^3\o} \\
	&		& \,\,\,\, Y^3\t\o h\t\o g\o\\
	&		& \,\,\,\, \tens w^3\underline{x^3\t} Y^3\t\t h\t\t g\t\\
	&	=	& Y^1aS(x^1Y^2)\alpha x^2\underline{X^1Y^3\o\o h\o\o} b S(w^1x^3\o \underline{X^2Y^3\o\t h\o\t})\\
	&		& \,\,\,\, \alpha w^2x^3\t\o \underline{X^3\o Y^3\t\o h\t\o} g\o\\
	&		& \,\,\,\, \tens w^3x^3\t\t \underline{X^3\t Y^3\t\t h\t\t} g\t\\
	&	=	& Y^1aS(x^1Y^2)\alpha x^2Y^3\o h\o X^1 b S(\underline{w^1x^3\o Y^3\t\o h\t\o} X^2)\\
	&		& \,\,\,\, \alpha \underline{w^2x^3\t\o Y^3\t\t\o h\t\t\o} X^3\o g\o\\
	&		& \,\,\,\, \tens \underline{w^3x^3\t\t Y^3\t\t\t h\t\t\t} X^3\t g\t\\	
	&	=	& \underline{Y^1}aS(\underline{x^1Y^2})\alpha \underline{x^2Y^3\o} h\o X^1 b S(w^1X^2) \alpha w^2 X^3\o g\o\\
	&		& \,\,\,\, \tens \underline{x^3Y^3\t} h^3\t w^3 X^3\t g\t\\
	&	=	& Y^1x^1\o aS(Y^2x^1\t)\alpha Y^3x^2 h\o \underline{X^1} b S(\underline{w^1X^2}) \alpha \underline{w^2X^3\o} g\o\\
	&		& \,\,\,\, \tens x^3 h^3\t \underline{w^3X^3\t} g\t\\		
	&	=	& Y^1x^1\o a S(Y^2x^1\t)\alpha Y^3x^2h\o X^1y^1\o bS(X^2y^1\t)\alpha X^3y^2 g\o\\
	&		& \,\,\,\, \tens x^3h\t y^3g\t\\
	&	=	& (q^1(x^1\la a)S(q^2)x^2h\o)(Q^1(y^1\la b)S(Q^2)y^2g\o) \tens (x^3h\t)(y^3g\t)\\
	&	=	& \chi(a\tens h)\chi(b\tens g)
\end{eqnarray*}

Next we show that $\chi$ is a coalgebra morphism.

\begin{eqnarray*}
\lefteqn{(\chi\tens\chi)\Delta(\chi^{-1}(b\tens 1))}\\
	&	=	& q^1(a^1w^1Y^1\la t^1T^1x^1\o b\o X^1\o \beta\o S(x^2X^2)\o g^1S(t^2R'\bt h^3T^3\t))S(q^2)\\
	&		&	\,\,\,\,\, a^2w^2\o W^1\o \underline{R\bt\o} y^2\o Y^3\o\o x^3\o\o X^3\o\o\\
	&		& \tens a^3w^2\t W^1\t \underline{R\bt\t} y^2\t Y^3\o\t x^3\o\t X^3\o\t\\
	&		& \tens Q^1(d^1w^3\o W^2\underline{R\bo} y^1Y^2t^3R'\bo \la h^1T^2x^1\t b\t \beta\t S(x^2X^2)\t g^2S(h^2T^3\o))S(Q^2)\\
	&		& \,\,\,\,\, d^2w^3\t\o W^3\o y^3\o Y^3\t\o x^3\t\o X^3\t\o\\
	&		& \tens d^3w^3\t\t W^3\t y^3\t Y^3\t\t x^3\t\t X^3\t\t\\
	&	=	& q^1(a^1w^1Y^1\la t^1T^1x^1\o b\o X^1\o \delta^1 S(t^2R'\bt h^3T^3\t x^2\t X^2\t))S(q^2)\\
	&		&	\,\,\,\,\, a^2w^2\o W^1\o u^1H^1R''\bt \underline{v^2 y^2\o} Y^3\o\o x^3\o\o X^3\o\o\\
	&		& \tens a^3w^2\t W^1\t u^2R\bt H^3\underline{v^3y^2\t} Y^3\o\t x^3\o\t X^3\o\t\\
	&		& \tens Q^1(d^1w^3\o W^2u^3R\bo H^2R''\bo \underline{v^1y^1}Y^2t^3R'\bo \la h^1T^2x^1\t b\t X^1\t\delta^2\\
	&		& \,\,\,\,\,  S(h^2T^3\o x^2\o X^t\o))S(Q^2)d^2w^3\t\o W^3\o \underline{y^3\o} Y^3\t\o x^3\t\o X^3\t\o\\
	&		& \tens d^3w^3\t\t W^3\t \underline{y^3\t} Y^3\t\t x^3\t\t X^3\t\t\\
	&	=	& q^1(a^1w^1Y^1\la t^1T^1x^1\o b\o X^1\o \delta^1 S(t^2R'\bt h^3T^3\t x^2\t X^2\t))S(q^2)\\
	&		&	\,\,\,\,\, a^2w^2\o W^1\o u^1\underline{H^1R''\bt y^1\t} v^2\underline{D^1Y^3\o\o x^3\o\o X^3\o\o}\\
	&		& \tens a^3w^2\t W^1\t u^2R\bt \underline{H^3y^2}v^3\o \underline{D^2Y^3\o\t x^3\o\t X^3\o\t}\\
	&		& \tens Q^1(d^1w^3\o W^2u^3R\bo \underline{H^2R''\bo y^1\o} v^1Y^2Y^2t^3R'\bo \la h^1T^2x^1\t b\t X^1\t\delta^2 \\
	&		& \,\,\,\,\,S(h^2T^3\o x^2\o X^t\o))S(Q^2)d^2w^3\t\o W^3\o \underline{y^3\o} v^3\t\o \underline{D^3\o Y^3\t\o x^3\t\o X^3\t\o}\\
	&		& \tens d^3w^3\t\t W^3\t \underline{y^3\t} v^3\t\t \underline{D^3\t Y^3\t\t x^3\t\t X^3\t\t}\\
	&	=	& q^1(a^1\o w^1\o \underline{Y^1\o} t^1T^1x^1\o b\o X^1\o \delta^1 S(a^1\t w^1\t \underline{Y^1\t} t^2R'\bt h^3T^3\t x^2\t X^2\t))\\
	&		& \,\,\,\,\, S(q^2)a^2w^2\o W^1\o u^1y^1H^1R''\bt \underline{v^2Y^3\o} x^3\o X^3\o D^1\\
	&		& \tens a^3w^2\t W^1\t u^2\underline{R\bt y^2\t} z^2H^3\o \underline{v^3\o Y^3\t\o} x^3\t\o X^3\t\o D^2\\
	&		& \tens Q^1(d^1w^3\o W^2u^3\underline{R\bo y^2\o} z^1H^2R''\bo \underline{v^1Y^2}t^3R'\bo \la h^1T^2x^1\t b\t X^1\t\delta^2 \\
	&		& \,\,\,\,\, S(h^2T^3\o x^2\o X^t\o))S(Q^2)d^2w^3\t\o W^3\o y^3\o z^3\o H^3\t\o \underline{v^3\t\o Y^3\t\t\o} \\
	&		& \,\,\,\,\, x^3\t\t\o X^3\t\t\o D^3\o\\
	&		& \tens d^3w^3\t\t W^3\t y^3\t z^3\t H^3\t\t \underline{v^3\t\t Y^3\t\t\t} x^3\t\t\t X^3\t\t\t D^3\t\\
	&	=	& q^1(a^1\o w^1\o Y^1\o F^1\o \underline{v^1\o\o t^1}T^1x^1\o b\o X^1\o \delta^1 \\
	&		& \,\,\,\,\, S(a^1\t w^1\t Y^1\t F^1\t \underline{v^1\o\t t^2R'\bt} h^3T^3\t x^2\t X^2\t))S(q^2)\\
	&		&	\,\,\,\,\, a^2w^2\o W^1\o u^1y^1H^1\underline{R''\bt Y^2\t} F^3v^2x^3\o X^3\o D^1\\
	&		& \tens a^3w^2\t W^1\t u^2y^2\o R\bt z^2H^3\o Y^3\o v^3\o x^3\t\o X^3\t\o D^2\\
	&		& \tens Q^1(d^1w^3\o W^2u^3y^2\t R\bo z^1H^2\underline{R''\bo Y^2\o} F^2 \underline{v^1\t t^3R'\bo} \la h^1T^2x^1\t b\t X^1\t\delta^2 \\
	&		& \,\,\,\,\, S(h^2T^3\o x^2\o X^t\o))S(Q^2)d^2w^3\t\o W^3\o y^3\o z^3\o H^3\t\o \\
	&		& \,\,\,\,\, Y^3\t\o v^3\t\o x^3\t\t\o X^3\t\t\o D^3\o\\
	&		& \tens d^3w^3\t\t W^3\t y^3\t z^3\t H^3\t\t Y^3\t\t v^3\t\t x^3\t\t\t X^3\t\t\t D^3\t\\
	&	=	& q^1(a^1\o w^1\o Y^1\o \underline{F^1\o t^1}v^1\o T^1x^1\o b\o X^1\o \delta^1 \\
	&		& \,\,\,\,\, S(a^1\t w^1\t Y^1\t \underline{F^1\t t^2}R'\bt v^1\t\t h^3T^3\t x^2\t X^2\t))S(q^2)\\
	&		&	\,\,\,\,\, a^2w^2\o W^1\o u^1y^1H^1Y^2\o R''\bt \underline{F^3}v^2x^3\o X^3\o D^1\\
	&		& \tens a^3w^2\t W^1\t u^2y^2\o R\bt z^2H^3\o Y^3\o v^3\o x^3\t\o X^3\t\o D^2\\
	&		& \tens Q^1(d^1w^3\o W^2u^3y^2\t R\bo z^1H^2Y^2\t R''\bo \underline{F^2 t^3}R'\bo v^1\t\o \la h^1T^2x^1\t b\t X^1\t\delta^2 \\
	&		& \,\,\,\,\, S(h^2T^3\o x^2\o X^t\o))S(Q^2)d^2w^3\t\o W^3\o y^3\o z^3\o H^3\t\o \\
	&		& \,\,\,\,\, Y^3\t\o v^3\t\o x^3\t\t\o X^3\t\t\o D^3\o\\
	&		& \tens d^3w^3\t\t W^3\t y^3\t z^3\t H^3\t\t Y^3\t\t v^3\t\t x^3\t\t\t X^3\t\t\t D^3\t\\
	&	=	& q^1(a^1\o w^1\o Y^1\o t^1U^1 v^1\o T^1x^1\o b\o X^1\o \delta^1 \\
	&		& \,\,\,\,\, S(a^1\t w^1\t Y^1\t t^2F^1\underline{U^2\o R'\bt} v^1\t\t h^3T^3\t x^2\t X^2\t))S(q^2)\\
	&		&	\,\,\,\,\, a^2w^2\o W^1\o u^1y^1H^1Y^2\o \underline{R''\bt t^3\t} F^3U^3v^2x^3\o X^3\o D^1\\
	&		& \tens a^3w^2\t W^1\t u^2y^2\o R\bt z^2H^3\o Y^3\o v^3\o x^3\t\o X^3\t\o D^2\\
	&		& \tens Q^1(d^1w^3\o W^2u^3y^2\t R\bo z^1H^2Y^2\t \underline{R''\bo t^3\o} F^2\underline{U^2\t R'\bo} v^1\t\o \la h^1T^2x^1\t b\t X^1\t\delta^2 \\
	&		& \,\,\,\,\, S(h^2T^3\o x^2\o X^t\o))S(Q^2)d^2w^3\t\o W^3\o y^3\o z^3\o H^3\t\o \\
	&		& \,\,\,\,\, Y^3\t\o v^3\t\o x^3\t\t\o X^3\t\t\o D^3\o\\
	&		& \tens d^3w^3\t\t W^3\t y^3\t z^3\t H^3\t\t Y^3\t\t v^3\t\t x^3\t\t\t X^3\t\t\t D^3\t\\
	&	=	& q^1(a^1\o w^1\o Y^1\o t^1U^1 v^1\o T^1x^1\o b\o X^1\o \delta^1 \\
	&		& \,\,\,\,\, S(a^1\t w^1\t Y^1\t t^2\underline{F^1R'\bt} U^2\t v^1\t\t h^3T^3\t x^2\t X^2\t))S(q^2)\\
	&		&	\,\,\,\,\, a^2w^2\o W^1\o u^1y^1H^1Y^2\o t^3\o \underline{R''\bt F^3}U^3v^2x^3\o X^3\o D^1\\
	&		& \tens a^3w^2\t W^1\t u^2y^2\o R\bt z^2H^3\o Y^3\o v^3\o x^3\t\o X^3\t\o D^2\\
	&		& \tens Q^1(d^1w^3\o W^2u^3y^2\t R\bo z^1H^2Y^2\t t^3\t \underline{R''\bo F^2R'\bo} U^2\o v^1\t\o \la h^1T^2x^1\t b\t X^1\t\delta^2 \\
	&		& \,\,\,\,\, S(h^2T^3\o x^2\o X^t\o))S(Q^2)d^2w^3\t\o W^3\o y^3\o z^3\o H^3\t\o \\
	&		& \,\,\,\,\, Y^3\t\o v^3\t\o x^3\t\t\o X^3\t\t\o D^3\o\\
	&		& \tens d^3w^3\t\t W^3\t y^3\t z^3\t H^3\t\t Y^3\t\t v^3\t\t x^3\t\t\t X^3\t\t\t D^3\t\\
	&	=	& q^1(\underline{a^1\o w^1\o} Y^1\o t^1U^1 v^1\o T^1x^1\o b\o X^1\o \delta^1 \\
	&		& \,\,\,\,\, S(\underline{a^1\t w^1\t} Y^1\t t^2F^1R'\bt\o A^2U^2\t v^1\t\t h^3T^3\t x^2\t X^2\t))S(q^2)\\
	&		&	\,\,\,\,\, \underline{a^2w^2\o} W^1\o u^1y^1H^1Y^2\o t^3\o F^2R'\bt\t A^3U^3v^2x^3\o X^3\o D^1\\
	&		& \tens \underline{a^3w^2\t} W^1\t u^2y^2\o R\bt z^2H^3\o Y^3\o v^3\o x^3\t\o X^3\t\o D^2\\
	&		& \tens Q^1(d^1\underline{w^3\o} W^2u^3y^2\t R\bo z^1H^2Y^2\t t^3\t F^3R'\bo A^1U^2\o v^1\t\o \la h^1T^2x^1\t b\t X^1\t\delta^2 \\
	&		& \,\,\,\,\, S(h^2T^3\o x^2\o X^t\o))S(Q^2)\\
	&		& \,\,\,\,\, d^2\underline{w^3\t\o} W^3\o y^3\o z^3\o H^3\t\o Y^3\t\o v^3\t\o x^3\t\t\o X^3\t\t\o D^3\o\\
	&		& \tens d^3\underline{w^3\t\t} W^3\t y^3\t z^3\t H^3\t\t Y^3\t\t v^3\t\t x^3\t\t\t X^3\t\t\t D^3\t\\
	&	=	& \underline{B^1a^1\o\o} w^1\o Y^1\o t^1U^1 v^1\o T^1x^1\o b\o X^1\o \delta^1 \\
	&		& \,\,\,\,\, S(\underline{B^2a^1\o\t} w^1\t Y^1\t t^2F^1R'\bt\o A^2U^2\t v^1\t\t h^3T^3\t x^2\t X^2\t)\alpha\\
	&		&	\,\,\,\,\, \underline{B^3a^1\t} w^2\underline{G^1W^1\o u^1y^1}H^1Y^2\o t^3\o F^2R'\bt\t A^3U^3v^2x^3\o X^3\o D^1\\
	&		& \tens a^2w^3\o \underline{G^2W^1\t u^2y^2\o} R\bt z^2H^3\o Y^3\o v^3\o x^3\t\o X^3\t\o D^2\\
	&		& \tens Q^1(d^1a^3\o w^3\t\o \underline{G^3\o W^2u^3y^2\t} R\bo z^1H^2Y^2\t t^3\t F^3R'\bo A^1U^2\o v^1\t\o \la h^1\\
	&		& \,\,\,\,\, T^2x^1\t b\t X^1\t\delta^2 S(h^2T^3\o x^2\o X^t\o))S(Q^2)d^2a^3\t\o w^3\t\t\o \\
	&		& \,\,\,\,\, \underline{G^3\t\o W^3\o y^3\o} z^3\o H^3\t\o Y^3\t\o v^3\t\o x^3\t\t\o X^3\t\t\o D^3\o\\
	&		& \tens d^3a^3\t\t w^3\t\t\t \underline{G^3\t\t W^3\t y^3\t} z^3\t H^3\t\t Y^3\t\t v^3\t\t x^3\t\t\t X^3\t\t\t D^3\t\\
	&	=	& a^1B^1\underline{w^1\o Y^1\o} t^1U^1 v^1\o T^1x^1\o b\o X^1\o \delta^1 \\
	&		& \,\,\,\,\, S(B^2\underline{w^1\t Y^1\t} t^2F^1R'\bt\o A^2U^2\t v^1\t\t h^3T^3\t x^2\t X^2\t)\alpha\\
	&		&	\,\,\,\,\, B^3\underline{w^2H^1Y^2\o} t^3\o F^2R'\bt\t A^3U^3v^2x^3\o X^3\o D^1\\
	&		& \tens a^2\underline{w^3\o G^1R\bt z^2H^3\o Y^3\o} v^3\o x^3\t\o X^3\t\o D^2\\
	&		& \tens Q^1(d^1a^3\o \underline{w^3\t\o G^2 R\bo z^1H^2Y^2\t} t^3\t F^3R'\bo A^1U^2\o v^1\t\o \la h^1T^2x^1\t b\t X^1\t\delta^2 \\
	&		& \,\,\,\,\, S(h^2T^3\o x^2\o X^t\o))S(Q^2)d^2a^3\t\o \\
	&		& \,\,\,\,\, \underline{w^3\t\t\o G^3\o z^3\o H^3\t\o Y^3\t\o} v^3\t\o x^3\t\t\o X^3\t\t\o D^3\o\\
	&		& \tens d^3a^3\t\t \underline{w^3\t\t\t G^3\t z^3\t H^3\t\t Y^3\t\t} v^3\t\t x^3\t\t\t X^3\t\t\t D^3\t\\
	&	=	& a^1\underline{B^1Y^1\o\o} w^1\o t^1U^1 v^1\o T^1x^1\o b\o X^1\o \delta^1 \\
	&		& \,\,\,\,\, S(\underline{B^2Y^1\o\t} w^1\t t^2F^1R'\bt\o A^2U^2\t v^1\t\t h^3T^3\t x^2\t X^2\t)\alpha\\
	&		&	\,\,\,\,\, \underline{B^3Y^1\t} w^2t^3\o F^2R'\bt\t A^3U^3v^2x^3\o X^3\o D^1\\
	&		& \tens a^2G^1R\bt z^2Y^3\o v^3\o x^3\t\o X^3\t\o D^2\\
	&		& \tens Q^1(d^1a^3\o G^2 R\bo z^1Y^2w^3 t^3\t F^3R'\bo A^1U^2\o v^1\t\o \la h^1T^2x^1\t b\t X^1\t\delta^2\\
	&		& \,\,\,\,\, S(h^2T^3\o x^2\o X^t\o))S(Q^2)d^2a^3\t\o \\
	&		& \,\,\,\,\, G^3\o z^3\o Y^3\t\o v^3\t\o x^3\t\t\o X^3\t\t\o D^3\o\\
	&		& \tens d^3a^3\t\t G^3\t z^3\t Y^3\t\t v^3\t\t x^3\t\t\t X^3\t\t\t D^3\t\\
	&	=	& a^1Y^1\underline{B^1w^1\o t^1}U^1 v^1\o T^1x^1\o b\o X^1\o \delta^1 \\
	&		& \,\,\,\,\, S(\underline{B^2 w^1\t t^2F^1}R'\bt\o A^2U^2\t v^1\t\t h^3T^3\t x^2\t X^2\t)\alpha\\
	&		&	\,\,\,\,\, \underline{B^3w^2t^3\o F^2}R'\bt\t A^3U^3v^2x^3\o X^3\o D^1\\
	&		& \tens a^2G^1R\bt z^2Y^3\o v^3\o x^3\t\o X^3\t\o D^2\\
	&		& \tens Q^1(d^1a^3\o G^2 R\bo z^1Y^2\underline{w^3 t^3\t F^3}R'\bo A^1U^2\o v^1\t\o \la h^1T^2x^1\t b\t X^1\t\delta^2 \\
	&		& \,\,\,\,\, S(h^2T^3\o x^2\o X^t\o))S(Q^2)d^2a^3\t\o G^3\o z^3\o \\
	&		& \,\,\,\,\, Y^3\t\o v^3\t\o x^3\t\t\o X^3\t\t\o D^3\o\\
	&		& \tens d^3a^3\t\t G^3\t z^3\t Y^3\t\t v^3\t\t x^3\t\t\t X^3\t\t\t D^3\t\\
	&	=	& a^1Y^1t^1U^1 v^1\o T^1x^1\o b\o X^1\o \delta^1 S(t^2\o R'\bt\o A^2U^2\t v^1\t\t h^3T^3\t x^2\t X^2\t)\alpha\\
	&		&	\,\,\,\,\, t^2\t R'\bt\t A^3U^3v^2x^3\o X^3\o D^1\\
	&		& \tens a^2G^1R\bt z^2Y^3\o v^3\o x^3\t\o X^3\t\o D^2\\
	&		& \tens Q^1(d^1a^3\o G^2 R\bo z^1Y^2t^3R'\bo A^1U^2\o v^1\t\o \la h^1T^2x^1\t b\t X^1\t\delta^2 \\
	&		& \,\,\,\,\, S(h^2T^3\o x^2\o X^t\o))S(Q^2)d^2a^3\t\o G^3\o z^3\o Y^3\t\o v^3\t\o x^3\t\t\o X^3\t\t\o D^3\o\\
	&		& \tens d^3a^3\t\t G^3\t z^3\t Y^3\t\t v^3\t\t x^3\t\t\t X^3\t\t\t D^3\t\\
	&	=	& a^1Y^1U^1 v^1\o T^1x^1\o b\o X^1\o \delta^1 S(A^2U^2\t v^1\t\t h^3T^3\t x^2\t X^2\t)\alpha\\
	&		&	\,\,\,\,\, A^3U^3v^2x^3\o X^3\o D^1\\
	&		& \tens a^2G^1R\bt z^2Y^3\o v^3\o x^3\t\o X^3\t\o D^2\\
	&		& \tens Q^1(\underline{d^1a^3\o} G^2 R\bo z^1Y^2A^1U^2\o v^1\t\o \la h^1T^2x^1\t b\t X^1\t\delta^2 S(h^2T^3\o x^2\o X^t\o))S(Q^2)\\
	&		& \,\,\,\,\, \underline{d^2a^3\t\o} G^3\o z^3\o Y^3\t\o v^3\t\o x^3\t\t\o X^3\t\t\o D^3\o\\
	&		& \tens \underline{d^3a^3\t\t} G^3\t z^3\t Y^3\t\t v^3\t\t x^3\t\t\t X^3\t\t\t D^3\t\\
	&	=	& a^1Y^1U^1 v^1\o T^1x^1\o b\o X^1\o \delta^1 S(A^2U^2\t v^1\t\t h^3T^3\t x^2\t X^2\t)\alpha\\
	&		&	\,\,\,\,\, A^3U^3v^2x^3\o X^3\o D^1\\
	&		& \tens a^2G^1R\bt z^2Y^3\o v^3\o x^3\t\o X^3\t\o D^2\\
	&		& \tens a^3\o Q^1(d^1 G^2 R\bo z^1Y^2A^1U^2\o v^1\t\o \la h^1T^2x^1\t b\t X^1\t\delta^2 S(h^2T^3\o x^2\o X^2\o))S(Q^2)\\
	&		& \,\,\,\,\, d^2G^3\o z^3\o Y^3\t\o v^3\t\o x^3\t\t\o X^3\t\t\o D^3\o\\
	&		& \tens a^3\t d^3G^3\t z^3\t Y^3\t\t v^3\t\t x^3\t\t\t X^3\t\t\t D^3\t\\
	&	=	& a^1Y^1U^1 \underline{v^1\o T^1x^1\o} b\o X^1\o \delta^1 S(A^2U^2\t \underline{v^1\t\t h^3T^3\t x^2\t} X^2\t)\alpha\\
	&		&	\,\,\,\,\, A^3U^3\underline{v^2x^3\o} X^3\o D^1\\
	&		& \tens a^2G^1R\bt z^2Y^3\o \underline{v^3\o x^3\t\o} X^3\t\o D^2\\
	&		& \tens a^3\o H^1d^1\o G^2\o R\bo\o z^1\o Y^2\o A^1\o U^2\o\o \underline{v^1\t\o\o h^1T^2x^1\t} b\t X^1\t\delta^2 \\
	&		& \,\,\,\,\, S(H^2d^1\t G^2\t R\bo\t z^1\t Y^2\t A^1\t U^2\o\t \underline{v^1\t\o\t h^2T^3\o x^2\o} X^2\o)\alpha\\
	&		& \,\,\,\,\, H^3d^2G^3\o z^3\o Y^3\t\o \underline{v^3\t\o x^3\t\t\o} X^3\t\t\o D^3\o\\
	&		& \tens a^3\t d^3G^3\t z^3\t Y^3\t\t \underline{v^3\t\t x^3\t\t\t} X^3\t\t\t D^3\t\\
	&	=	& a^1Y^1U^1 T^1v^1\o x^1\o b\o X^1\o \delta^1 S(A^2\underline{U^2\t h^3}T^3\t v^2\t x^2\o\t w^1\t X^2\t)\alpha\\
	&		&	\,\,\,\,\, A^3U^3v^3x^2\t w^2 X^3\o D^1\\
	&		& \tens a^2G^1R\bt z^2Y^3\o x^3\o w^3\o X^3\t\o D^2\\
	&		& \tens a^3\o H^1d^1\o G^2\o R\bo\o z^1\o Y^2\o A^1\o \underline{U^2\o\o h^1}T^2v^1\t x^1\t b\t X^1\t\delta^2 \\
	&		& \,\,\,\,\, S(H^2d^1\t G^2\t R\bo\t z^1\t Y^2\t A^1\t \underline{U^2\o\t h^2}T^3\o v^2\o x^2\o\o w^1\o X^2\o)\alpha\\
	&		& \,\,\,\,\, H^3d^2G^3\o z^3\o Y^3\t\o x^3\t\o w^3\t\o X^3\t\t\o D^3\o\\
	&		& \tens a^3\t d^3G^3\t z^3\t Y^3\t\t x^3\t\t w^3\t\t X^3\t\t\t D^3\t\\
	&	=	& a^1Y^1U^1 T^1v^1\o x^1\o b\o X^1\o \delta^1 S(\underline{A^2h^3} U^2\t\t T^3\t v^2\t x^2\o\t w^1\t X^2\t)\alpha\\
	&		&	\,\,\,\,\, \underline{A^3}U^3v^3x^2\t w^2 X^3\o D^1\\
	&		& \tens a^2G^1R\bt z^2Y^3\o x^3\o w^3\o X^3\t\o D^2\\
	&		& \tens a^3\o H^1d^1\o G^2\o R\bo\o z^1\o Y^2\o \underline{A^1\o h^1} U^2\o T^2v^1\t x^1\t b\t X^1\t\delta^2 \\
	&		& \,\,\,\,\, S(H^2d^1\t G^2\t R\bo\t z^1\t Y^2\t \underline{A^1\t h^2}U^2\t\o T^3\o v^2\o x^2\o\o w^1\o X^2\o)\alpha\\
	&		& \,\,\,\,\, H^3d^2G^3\o z^3\o Y^3\t\o x^3\t\o w^3\t\o X^3\t\t\o D^3\o\\
	&		& \tens a^3\t d^3G^3\t z^3\t Y^3\t\t x^3\t\t w^3\t\t X^3\t\t\t D^3\t\\
	&	=	& a^1Y^1\underline{U^1 T^1v^1\o} x^1\o b\o X^1\o \delta^1 S(B^2\underline{A^2\t U^2\t\t T^3\t v^2\t} x^2\o\t w^1\t X^2\t)\alpha \\
	&		& \,\,\,\,\, B^3\underline{A^3U^3v^3}x^2\t w^2 X^3\o D^1\\
	&		& \tens a^2G^1R\bt z^2Y^3\o x^3\o w^3\o X^3\t\o D^2\\
	&		& \tens a^3\o H^1d^1\o G^2\o R\bo\o z^1\o Y^2\o \underline{A^1 U^2\o T^2v^1\t} x^1\t b\t X^1\t\delta^2 \\
	&		& \,\,\,\,\, S(H^2d^1\t G^2\t R\bo\t z^1\t Y^2\t B^1\underline{A^2\o U^2\t\o T^3\o v^2\o} x^2\o\o w^1\o X^2\o)\alpha\\
	&		& \,\,\,\,\, H^3d^2G^3\o z^3\o Y^3\t\o x^3\t\o w^3\t\o X^3\t\t\o D^3\o\\
	&		& \tens a^3\t d^3G^3\t z^3\t Y^3\t\t x^3\t\t w^3\t\t X^3\t\t\t D^3\t\\
	&	=	& a^1Y^1U^1x^1\o b\o X^1\o \delta^1 S(\underline{B^2U^3\o\t x^2\o\t} w^1\t X^2\t)\alpha \underline{B^3U^3\t x^2\t} w^2 X^3\o D^1\\
	&		& \tens a^2G^1R\bt z^2Y^3\o x^3\o w^3\o X^3\t\o D^2\\
	&		& \tens a^3\o H^1d^1\o G^2\o R\bo\o z^1\o Y^2\o U^2x^1\t b\t \delta^2 S(H^2d^1\t G^2\t R\bo\t z^1\t Y^2\t\\
	&		& \,\,\,\,\,  \underline{B^1U^3\o\o x^2\o\o} w^1\o X^2\o)\alpha\\
	&		& \,\,\,\,\, H^3d^2G^3\o z^3\o Y^3\t\o x^3\t\o w^3\t\o X^3\t\t\o D^3\o\\
	&		& \tens a^3\t d^3G^3\t z^3\t Y^3\t\t x^3\t\t w^3\t\t X^3\t\t\t D^3\t\\
	&	=	& a^1\underline{Y^1U^1x^1\o} b\o X^1\o \delta^1 S(B^2 w^1\t X^2\t)\alpha B^3w^2 X^3\o D^1\\
	&		& \tens a^2G^1R\bt z^2\underline{Y^3\o x^3\o} w^3\o X^3\t\o D^2\\
	&		& \tens a^3\o H^1d^1\o G^2\o R\bo\o z^1\o \underline{Y^2\o U^2x^1\t} b\t X^1\t\delta^2 S(H^2d^1\t G^2\t R\bo\t z^1\t \\
	&		& \,\,\,\,\, \underline{Y^2\t U^3x^2}B^1w^1\o X^2\o)\alpha\\
	&		& \,\,\,\,\, H^3d^2G^3\o z^3\o \underline{Y^3\t\o x^3\t\o} w^3\t\o X^3\t\t\o D^3\o\\
	&		& \tens a^3\t d^3G^3\t z^3\t \underline{Y^3\t\t x^3\t\t} w^3\t\t X^3\t\t\t D^3\t\\
	&	=	& a^1Y^1 b\o X^1\o \delta^1 S(B^2 w^1\t X^2\t)\alpha B^3w^2 X^3\o D^1\\
	&		& \tens a^2\underline{G^1}R\bt z^2x^3\o Y^3\t\o w^3\o X^3\t\o D^2\\
	&		& \tens a^3\o H^1\underline{d^1\o G^2\o} R\bo\o z^1\o x^1Y^2 b\t \delta^2 S(H^2\underline{d^1\t G^2\t} R\bo\t z^1\t x^2Y^3\o B^1w^1\o X^2\o)\alpha\\
	&		& \,\,\,\,\, H^3\underline{d^2G^3\o} z^3\o x^3\t\o Y^3\t\t\o w^3\t\o X^3\t\t\o D^3\o\\
	&		& \tens a^3\t \underline{d^3G^3\t} z^3\t x^3\t\t Y^3\t\t\t w^3\t\t X^3\t\t\t D^3\t\\
	&	=	& a^1Y^1 b\o X^1\o \delta^1 S(B^2 w^1\t X^2\t)\alpha B^3w^2 X^3\o D^1\\
	&		& \tens a^2G^1A^1\underline{d^1\o R\bt} z^2x^3\o Y^3\t\o w^3\o X^3\t\o D^2\\
	&		& \tens a^3\o H^1G^2\o\o A^2\o \underline{d^1\t\o R\bo\o} z^1\o x^1Y^2 b\t X^1\t\delta^2 S(H^2G^2\o\t A^2\t \\
	&		& \,\,\,\,\, \underline{d^1\t\t R\bo\t} z^1\t x^2Y^3\o B^1w^1\o X^2\o)\alpha\\
	&		& \,\,\,\,\, H^3G^2\t A^3 d^2 z^3\o x^3\t\o Y^3\t\t\o w^3\t\o X^3\t\t\o D^3\o\\
	&		& \tens a^3\t G^3d^3 z^3\t x^3\t\t Y^3\t\t\t w^3\t\t X^3\t\t\t D^3\t\\
	&	=	& a^1Y^1 b\o X^1\o \delta^1 S(B^2 w^1\t X^2\t)\alpha B^3w^2 X^3\o D^1\\
	&		& \tens a^2G^1A^1R\bt d^1\t z^2x^3\o Y^3\t\o w^3\o X^3\t\o D^2\\
	&		& \tens a^3\o \underline{H^1G^2\o\o} A^2\o R\bo\o d^1\o\o z^1\o x^1Y^2 b\t X^1\t\delta^2 \\
	&		& \,\,\,\,\, S(\underline{H^2G^2\o\t} A^2\t R\bo\t d^1\o\t z^1\t x^2Y^3\o B^1w^1\o X^2\o)\alpha\\
	&		& \,\,\,\,\, \underline{H^3G^2\t} A^3 d^2 z^3\o x^3\t\o Y^3\t\t\o w^3\t\o X^3\t\t\o D^3\o\\
	&		& \tens a^3\t G^3d^3 z^3\t x^3\t\t Y^3\t\t\t w^3\t\t X^3\t\t\t D^3\t\\
	&	=	& a^1Y^1 b\o X^1\o \delta^1 S(B^2 w^1\t X^2\t)\alpha B^3w^2 X^3\o D^1\\
	&		& \tens a^2G^1\underline{A^1}R\bt d^1\t z^2x^3\o Y^3\t\o w^3\o X^3\t\o D^2\\
	&		& \tens a^3\o G^2\underline{H^1 A^2\o} R\bo\o d^1\o\o z^1\o x^1Y^2 b\t X^1\t\delta^2 \\
	&		& \,\,\,\,\, S(\underline{H^2 A^2\t} R\bo\t d^1\o\t z^1\t x^2Y^3\o B^1w^1\o X^2\o)\alpha\\
	&		& \,\,\,\,\, \underline{H^3A^3} d^2 z^3\o x^3\t\o Y^3\t\t\o w^3\t\o X^3\t\t\o D^3\o\\
	&		& \tens a^3\t G^3d^3 z^3\t x^3\t\t Y^3\t\t\t w^3\t\t X^3\t\t\t D^3\t\\
	&	=	& a^1Y^1 b\o X^1\o \delta^1 S(B^2 w^1\t X^2\t)\alpha B^3w^2 X^3\o D^1\\
	&		& \tens a^2G^1A^1\o y^1R\bt \underline{d^1\t z^2}x^3\o Y^3\t\o w^3\o X^3\t\o D^2\\
	&		& \tens a^3\o G^2A^1\t y^2 R\bo\o \underline{d^1\o\o z^1\o} x^1Y^2 b\t X^1\t\delta^2 \\
	&		& \,\,\,\,\, S(A^2y^3 R\bo\t \underline{d^1\o\t z^1\t} x^2Y^3\o B^1w^1\o X^2\o)\alpha\\
	&		& \,\,\,\,\, A^3 \underline{d^2 z^3\o} x^3\t\o Y^3\t\t\o w^3\t\o X^3\t\t\o D^3\o\\
	&		& \tens a^3\t G^3\underline{d^3 z^3\t} x^3\t\t Y^3\t\t\t w^3\t\t X^3\t\t\t D^3\t\\
	&	=	& a^1Y^1 b\o X^1\o \delta^1 S(B^2 w^1\t X^2\t)\alpha B^3w^2 X^3\o D^1\\
	&		& \tens a^2G^1A^1\o y^1R\bt d^2z^2\o \underline{h^1x^3\o Y^3\t\o w^3\o X^3\t\o }D^2\\
	&		& \tens a^3\o G^2A^1\t y^2 R\bo\o d^1\o z^1\o  x^1Y^2 b\t X^1\t\delta^2 \\
	&		& \,\,\,\,\, S(A^2y^3 R\bo\t d^1\t z^1\t x^2Y^3\o B^1w^1\o X^2\o)\alpha\\
	&		& \,\,\,\,\, A^3 d^3z^2\t \underline{h^2 x^3\t\o Y^3\t\t\o w^3\t\o X^3\t\t\o}D^3\o\\
	&		& \tens a^3\t G^3z^3\underline{h^3 x^3\t\t Y^3\t\t\t w^3\t\t X^3\t\t\t} D^3\t\\
	&	=	& a^1Y^1 b\o X^1\o \delta^1 S(B^2 w^1\t X^2\t)\alpha B^3w^2 X^3\o D^1\\
	&		& \tens a^2G^1A^1\o y^1R\bt d^2\underline{z^2\o x^3\o\o} Y^3\t\o\o w^3\o\o X^3\t\o\o h^1D^2\\
	&		& \tens a^3\o G^2A^1\t y^2 R\bo\o d^1\o \underline{z^1\o  x^1}Y^2 b\t X^1\t\delta^2 \\
	&		& \,\,\,\,\, S(A^2y^3 R\bo\t d^1\t \underline{z^1\t x^2}Y^3\o B^1w^1\o X^2\o)\alpha\\
	&		& \,\,\,\,\, A^3 d^3\underline{z^2\t x^3\o\t} Y^3\t\o\t w^3\o\t X^3\t\o\t h^2D^3\o\\
	&		& \tens a^3\t G^3\underline{z^3x^3\t} Y^3\t\t w^3\t X^3\t\t h^3 D^3\t\\
	&	=	& a^1Y^1 b\o X^1\o \delta^1 S(B^2 w^1\t X^2\t)\alpha B^3w^2 X^3\o D^1\\
	&		& \tens a^2G^1A^1\o y^1R\bt \underline{d^2z^3\o} x^2\t\o k^2\o Y^3\t\o\o w^3\o\o X^3\t\o\o h^1D^2\\
	&		& \tens a^3\o G^2A^1\t y^2 R\bo\o \underline{d^1\o z^1}x^1Y^2 b\t X^1\t\delta^2 \\
	&		& \,\,\,\,\, S(A^2y^3 R\bo\t \underline{d^1\t z^2}x^2\o k^1Y^3\o B^1w^1\o X^2\o)\alpha\\
	&		& \,\,\,\,\, A^3 \underline{d^3z^3\t} x^2\t\t k^2\t Y^3\t\o\t w^3\o\t X^3\t\o\t h^2D^3\o\\
	&		& \tens a^3\t G^3x^3k^3 Y^3\t\t w^3\t X^3\t\t h^3 D^3\t\\
	&	=	& a^1Y^1 b\o X^1\o \delta^1 S(B^2 w^1\t X^2\t)\alpha B^3w^2 X^3\o D^1\\
	&		& \tens a^2G^1A^1\o \underline{y^1R\bt d^3}z^2\t u^2 x^2\t\o k^2\o Y^3\t\o\o w^3\o\o X^3\t\o\o h^1D^2\\
	&		& \tens a^3\o G^2A^1\t \underline{y^2 R\bo\o d^1}z^1x^1Y^2 b\t X^1\t\delta^2 \\
	&		& \,\,\,\,\, S(A^2\underline{y^3 R\bo\t d^2}z^2\o u^1x^2\o k^1Y^3\o B^1w^1\o X^2\o)\alpha\\
	&		& \,\,\,\,\, A^3 z^3u^3 x^2\t\t k^2\t Y^3\t\o\t w^3\o\t X^3\t\o\t h^2D^3\o\\
	&		& \tens a^3\t G^3x^3k^3 Y^3\t\t w^3\t X^3\t\t h^3 D^3\t\\
	&	=	& a^1Y^1 b\o X^1\o \delta^1 S(B^2 w^1\t X^2\t)\alpha B^3w^2 X^3\o D^1\\
	&		& \tens a^2G^1\underline{A^1\o R\bt} y^2\underline{R'\bt z^2\t} u^2 x^2\t\o k^2\o Y^3\t\o\o w^3\o\o X^3\t\o\o h^1D^2\\
	&		& \tens a^3\o G^2\underline{A^1\t R\bo} y^1 z^1x^1Y^2 b\t X^1\t\delta^2 \\
	&		& \,\,\,\,\, S(A^2y^3\underline{R'\bo z^2\o} u^1x^2\o k^1Y^3\o B^1w^1\o X^2\o)\alpha\\
	&		& \,\,\,\,\, A^3 z^3u^3 x^2\t\t k^2\t Y^3\t\o\t w^3\o\t X^3\t\o\t h^2D^3\o\\
	&		& \tens a^3\t G^3x^3k^3 Y^3\t\t w^3\t X^3\t\t h^3 D^3\t\\
	&	=	& a^1Y^1 b\o X^1\o \delta^1 S(B^2 w^1\t X^2\t)\alpha B^3w^2 X^3\o D^1\\
	&		& \tens a^2G^1R\bt \underline{A^1\t y^2z^2\o} R'\bt u^2 x^2\t\o k^2\o Y^3\t\o\o w^3\o\o X^3\t\o\o h^1D^2\\
	&		& \tens a^3\o G^2R\bo \underline{A^1\o y^1 z^1}x^1Y^2 b\t X^1\t\delta^2 S(\underline{A^2y^3z^2\t} R'\bo u^1x^2\o k^1Y^3\o B^1w^1\o X^2\o)\alpha\\
	&		& \,\,\,\,\, \underline{A^3 z^3}u^3 x^2\t\t k^2\t Y^3\t\o\t w^3\o\t X^3\t\o\t h^2D^3\o\\
	&		& \tens a^3\t G^3x^3k^3 Y^3\t\t w^3\t X^3\t\t h^3 D^3\t\\
	&	=	& a^1Y^1 b\o X^1\o \delta^1 S(B^2 w^1\t X^2\t)\alpha B^3w^2 X^3\o D^1\\
	&		& \tens a^2G^1R\bt A^1 R'\bt u^2 x^2\t\o k^2\o Y^3\t\o\o w^3\o\o X^3\t\o\o h^1D^2\\
	&		& \tens a^3\o G^2R\bo x^1Y^2 b\t X^1\t\delta^2 S(A^2 R'\bo u^1x^2\o k^1Y^3\o B^1w^1\o X^2\o)\alpha\\
	&		& \,\,\,\,\, A^3u^3 x^2\t\t k^2\t Y^3\t\o\t w^3\o\t X^3\t\o\t h^2D^3\o\\
	&		& \tens a^3\t G^3x^3k^3 Y^3\t\t w^3\t X^3\t\t h^3 D^3\t\\
	&	=	& a^1Y^1 b\o \underline{X^1\o} \delta^1 S(B^2 \underline{w^1\t X^2\t})\alpha B^3\underline{w^2 X^3\o} D^1\\
	&		& \tens a^2G^1R\bt A^1 R'\bt u^2 x^2\t\o k^2\o Y^3\t\o\o \underline{w^3\o\o X^3\t\o\o} h^1D^2\\
	&		& \tens a^3\o G^2R\bo x^1Y^2 b\t \underline{X^1\t}\delta^2 S(t^1B^1\underline{w^1\o X^2\o})\alpha t^2\beta S(A^2 R'\bo u^1x^2\o k^1Y^3\o t^3)\alpha\\
	&		& \,\,\,\,\, A^3u^3 x^2\t\t k^2\t Y^3\t\o\t \underline{w^3\o\t X^3\t\o\t} h^2D^3\o\\
	&		& \tens a^3\t G^3x^3k^3 Y^3\t\t \underline{w^3\t X^3\t\t} h^3 D^3\t\\
	&	=	& a^1Y^1 b\o W^1\o X^1\o \delta^1 S(\underline{B^2 W^2\o\t} X^2\t)\alpha \underline{B^3W^2\t} X^3D^1\\
	&		& \tens a^2G^1R\bt A^1 R'\bt u^2 x^2\t\o k^2\o Y^3\t\o\o W^3\o\o h^1D^2\\
	&		& \tens a^3\o G^2R\bo x^1Y^2 b\t W^1\t X^1\t\delta^2 S(t^1\underline{B^1W^2\o\o} X^2\o)\alpha t^2\beta S(A^2 R'\bo u^1x^2\o k^1Y^3\o t^3)\alpha\\
	&		& \,\,\,\,\, A^3u^3 x^2\t\t k^2\t Y^3\t\o\t W^3\o\t h^2D^3\o\\
	&		& \tens a^3\t G^3x^3k^3 Y^3\t\t W^3\t h^3 D^3\t\\
	&	=	& a^1Y^1 b\o W^1\o X^1\o \delta^1 S(B^2 X^2\t)\alpha B^3 X^3D^1\\
	&		& \tens a^2G^1R\bt A^1 R'\bt u^2 x^2\t\o k^2\o Y^3\t\o\o W^3\o\o h^1D^2\\
	&		& \tens a^3\o G^2R\bo x^1Y^2 b\t W^1\t X^1\t\delta^2 S(\underline{t^1W^2}B^1X^2\o)\alpha \underline{t^2}\beta S(A^2 R'\bo u^1x^2\o k^1Y^3\o \underline{t^3})\alpha\\
	&		& \,\,\,\,\, A^3u^3 x^2\t\t k^2\t Y^3\t\o\t W^3\o\t h^2D^3\o\\
	&		& \tens a^3\t G^3x^3k^3 Y^3\t\t W^3\t h^3 D^3\t\\
	&	=	& a^1Y^1 b\o W^1\o V^1\o X^1\o \delta^1 S(B^2 X^2\t)\alpha B^3 X^3\underline{T^1}d^1D^1\\
	&		& \tens a^2G^1R\bt A^1 R'\bt u^2 x^2\t\o k^2\o Y^3\t\o\o W^3\o\o h^1D^2\\
	&		& \tens a^3\o G^2R\bo x^1Y^2 b\t W^1\t V^1\t X^1\t\delta^2 S(t^1V^2B^1X^2\o)\alpha t^2V^3\o T^2d^2\beta \\
	&		& \,\,\,\,\, S(A^2 R'\bo u^1x^2\o k^1Y^3\o W^2t^3V^3\t T^3d^3)\alpha A^3u^3 x^2\t\t k^2\t Y^3\t\o\t W^3\o\t h^2D^3\o\\
	&		& \tens a^3\t G^3x^3k^3 Y^3\t\t W^3\t h^3 D^3\t\\
	&	=	& a^1Y^1 b\o W^1\o V^1\o \underline{X^1\o T^1\o\o\o}\delta^1 S(B^2 \underline{X^2\t T^1\o\t\t})\alpha B^3 \underline{X^3T^1\t} d^1D^1\\
	&		& \tens a^2G^1R\bt A^1 R'\bt u^2 x^2\t\o k^2\o Y^3\t\o\o W^3\o\o h^1D^2\\
	&		& \tens a^3\o G^2R\bo x^1Y^2 b\t W^1\t V^1\t \underline{X^1\t T^1\o\o\t}\delta^2 S(t^1V^2B^1\underline{X^2\o T^1\o\t\o})\alpha t^2V^3\o T^2d^2\beta \\
	&		& \,\,\,\,\, S(A^2 R'\bo u^1x^2\o k^1Y^3\o W^2t^3V^3\t T^3d^3)\alpha A^3u^3 x^2\t\t k^2\t Y^3\t\o\t W^3\o\t h^2D^3\o\\
	&		& \tens a^3\t G^3x^3k^3 Y^3\t\t W^3\t h^3 D^3\t\\
	&	=	& a^1Y^1 b\o W^1\o V^1\o T^1\o\o X^1\o\delta^1 S(\underline{B^2 T^1\t\o\t} X^2\t )\alpha \underline{B^3 T^1\t\t} X^3 d^1D^1\\
	&		& \tens a^2G^1R\bt A^1 R'\bt u^2 x^2\t\o k^2\o Y^3\t\o\o W^3\o\o h^1D^2\\
	&		& \tens a^3\o G^2R\bo x^1Y^2 b\t W^1\t V^1\t T^1\o\t X^1\t\delta^2 S(t^1V^2\underline{B^1T^1\t\o\o} X^2\o )\alpha t^2V^3\o T^2d^2\beta \\
	&		& \,\,\,\,\, S(A^2 R'\bo u^1x^2\o k^1Y^3\o W^2t^3V^3\t T^3d^3)\alpha A^3u^3 x^2\t\t k^2\t Y^3\t\o\t W^3\o\t h^2D^3\o\\
	&		& \tens a^3\t G^3x^3k^3 Y^3\t\t W^3\t h^3 D^3\t\\
	&	=	& a^1Y^1 b\o W^1\o V^1\o T^1\o\o \underline{X^1\o}\delta^1 S(B^2 \underline{X^2\t} )\alpha B^3 \underline{X^3} d^1D^1\\
	&		& \tens a^2G^1R\bt A^1 R'\bt u^2 x^2\t\o k^2\o Y^3\t\o\o W^3\o\o h^1D^2\\
	&		& \tens a^3\o G^2R\bo x^1Y^2 b\t W^1\t V^1\t T^1\o\t \underline{X^1\t}\delta^2 S(t^1V^2T^1\t B^1 \underline{X^2\o} )\alpha t^2V^3\o T^2d^2\beta \\
	&		& \,\,\,\,\, S(A^2 R'\bo u^1x^2\o k^1Y^3\o W^2t^3V^3\t T^3d^3)\alpha A^3u^3 x^2\t\t k^2\t Y^3\t\o\t W^3\o\t h^2D^3\o\\
	&		& \tens a^3\t G^3x^3k^3 Y^3\t\t W^3\t h^3 D^3\t\\
	&	=	& a^1Y^1 b\o W^1\o V^1\o T^1\o\o z^1\o X^1\o H^1\o\o \delta^1 \\
	&		& \,\,\,\,\, S(\underline{B^2 z^2\o\t} w^1\t X^2\t H^1\t\t )\alpha \underline{B^3 z^2\t} w^2X^3\o H^2 d^1D^1\\
	&		& \tens a^2G^1R\bt A^1 R'\bt u^2 x^2\t\o k^2\o Y^3\t\o\o W^3\o\o h^1D^2\\
	&		& \tens a^3\o G^2R\bo x^1Y^2 b\t W^1\t V^1\t T^1\o\t z^1\t X^1\t H^1\o\t\delta^2 \\
	&		& \,\,\,\,\, S(t^1V^2T^1\t \underline{B^1 z^2\o\o} w^1\o X^2\o H^1\t\o )\alpha t^2V^3\o T^2z^3w^3X^3\t H^3d^2\beta \\
	&		& \,\,\,\,\, S(A^2 R'\bo u^1x^2\o k^1Y^3\o W^2t^3V^3\t T^3d^3)\alpha A^3u^3 x^2\t\t k^2\t Y^3\t\o\t W^3\o\t h^2D^3\o\\
	&		& \tens a^3\t G^3x^3k^3 Y^3\t\t W^3\t h^3 D^3\t\\
	&	=	& a^1Y^1 b\o W^1\o \underline{V^1\o T^1\o\o z^1\o} X^1\o  \delta^1 S(B^2 w^1\t X^2\t )\alpha B^3 w^2X^3\o d^1D^1\\
	&		& \tens a^2G^1R\bt A^1 R'\bt u^2 x^2\t\o k^2\o Y^3\t\o\o W^3\o\o h^1D^2\\
	&		& \tens a^3\o G^2R\bo x^1Y^2 b\t W^1\t \underline{V^1\t T^1\o\t z^1\t} X^1\t \delta^2 \\
	&		& \,\,\,\,\, S(\underline{t^1V^2T^1\t z^2}B^1 w^1\o X^2\o )\alpha \underline{t^2V^3\o T^2z^3}w^3X^3\t d^2\beta \\
	&		& \,\,\,\,\, S(A^2 R'\bo u^1x^2\o k^1Y^3\o W^2\underline{t^3V^3\t T^3}d^3)\alpha A^3u^3 x^2\t\t k^2\t Y^3\t\o\t W^3\o\t h^2D^3\o\\
	&		& \tens a^3\t G^3x^3k^3 Y^3\t\t W^3\t h^3 D^3\t\\
	&	=	& a^1Y^1 b\o W^1\o T^1\o X^1\o  \delta^1 S(B^2 w^1\t X^2\t )\alpha B^3 w^2X^3\o d^1D^1\\
	&		& \tens a^2G^1R\bt A^1 R'\bt u^2 x^2\t\o k^2\o Y^3\t\o\o W^3\o\o h^1D^2\\
	&		& \tens a^3\o G^2R\bo x^1Y^2 b\t W^1\t T^1\t X^1\t \delta^2 S(T^2\o B^1 w^1\o X^2\o )\alpha T^2\t w^3X^3\t d^2\beta \\
	&		& \,\,\,\,\, S(A^2 R'\bo u^1x^2\o k^1Y^3\o W^2T^3d^3)\alpha A^3u^3 x^2\t\t k^2\t Y^3\t\o\t W^3\o\t h^2D^3\o\\
	&		& \tens a^3\t G^3x^3k^3 Y^3\t\t W^3\t h^3 D^3\t\\
	&	=	& a^1Y^1 b\o W^1\o \underline{X^1\o  \delta^1 S(B^2 w^1\t X^2\t )\alpha B^3 w^2X^3\o} d^1D^1\\
	&		& \tens a^2G^1R\bt A^1 R'\bt u^2 x^2\t\o k^2\o Y^3\t\o\o W^3\o\o h^1D^2\\
	&		& \tens a^3\o G^2R\bo x^1Y^2 b\t W^1\t \underline{X^1\t \delta^2 S(B^1 w^1\o X^2\o )\alpha w^3X^3\t} d^2\beta \\
	&		& \,\,\,\,\, S(A^2 R'\bo u^1x^2\o k^1Y^3\o W^2d^3)\alpha A^3u^3 x^2\t\t k^2\t Y^3\t\o\t W^3\o\t h^2D^3\o\\
	&		& \tens a^3\t G^3x^3k^3 Y^3\t\t W^3\t h^3 D^3\t\\
	&	=	& a^1Y^1 b\o \underline{W^1\o} d^1D^1\\
	&		& \tens a^2G^1R\bt A^1 R'\bt u^2 x^2\t\o k^2\o Y^3\t\o\o \underline{W^3\o\o} h^1D^2\\
	&		& \tens a^3\o G^2R\bo x^1Y^2 b\t \underline{W^1\t} d^2\beta S(A^2 R'\bo u^1x^2\o k^1Y^3\o \underline{W^2}d^3)\alpha \\
	&		& \,\,\,\,\, A^3u^3 x^2\t\t k^2\t Y^3\t\o\t \underline{W^3\o\t} h^2D^3\o\\
	&		& \tens a^3\t G^3x^3k^3 Y^3\t\t \underline{W^3\t} h^3 D^3\t\\
	&	=	& a^1Y^1 b\o X^1\o T^1\o \underline{t^1\o\o d^1}D^1\\
	&		& \tens a^2G^1R\bt A^1 R'\bt u^2 x^2\t\o k^2\o Y^3\t\o\o W^2\o X^2\t\o T^3\o t^2\o h^1D^2\\
	&		& \tens a^3\o G^2R\bo x^1Y^2 b\t X^1\t T^1\t \underline{t^1\o\t d^2}\beta S(A^2 R'\bo u^1x^2\o k^1Y^3\o W^1X^2\o T^2\underline{t^1\t d^3})\alpha \\
	&		& \,\,\,\,\, A^3u^3 x^2\t\t k^2\t Y^3\t\o\t W^2\t X^2\t\t T^3\t t^2\t h^2D^3\o\\
	&		& \tens a^3\t G^3x^3k^3 Y^3\t\t W^3X^3t^3 h^3 D^3\t\\
	&	=	& a^1Y^1 b\o X^1\o \underline{T^1\o d^1}t^1D^1\\
	&		& \tens a^2G^1R\bt A^1 R'\bt u^2 x^2\t\o k^2\o Y^3\t\o\o W^2\o X^2\t\o \underline{T^3\o} t^2\o h^1D^2\\
	&		& \tens a^3\o G^2R\bo x^1Y^2 b\t X^1\t \underline{T^1\t d^2}\beta S(A^2 R'\bo u^1x^2\o k^1Y^3\o W^1X^2\o \underline{T^2d^3})\alpha \\
	&		& \,\,\,\,\, A^3u^3 x^2\t\t k^2\t Y^3\t\o\t W^2\t X^2\t\t \underline{T^3\t} t^2\t h^2D^3\o\\
	&		& \tens a^3\t G^3x^3k^3 Y^3\t\t W^3X^3t^3 h^3 D^3\t\\
	&	=	& a^1Y^1 b\o X^1\o d^1t^1D^1\\
	&		& \tens a^2G^1R\bt A^1 R'\bt u^2 x^2\t\o \underline{k^2\o Y^3\t\o\o W^2\o} X^2\t\o d^3\t\o T^3\o t^2\o h^1D^2\\
	&		& \tens a^3\o G^2R\bo x^1Y^2 b\t X^1\t d^2T^1\beta S(A^2 R'\bo u^1x^2\o \underline{k^1Y^3\o W^1}X^2\o d^3\o T^2)\alpha \\
	&		& \,\,\,\,\, A^3u^3 x^2\t\t \underline{k^2\t Y^3\t\o\t W^2\t} X^2\t\t d^3\t\t T^3\t t^2\t h^2D^3\o\\
	&		& \tens a^3\t G^3x^3\underline{k^3 Y^3\t\t W^3}X^3t^3 h^3 D^3\t\\
	&	=	& a^1Y^1 b\o \underline{X^1\o d^1}t^1D^1\\
	&		& \tens a^2G^1R\bt A^1 R'\bt u^2 x^2\t\o Y^3\o\t\o \underline{X^2\t\o d^3\t\o} T^3\o t^2\o h^1D^2\\
	&		& \tens a^3\o G^2R\bo x^1Y^2 b\t \underline{X^1\t d^2}T^1\beta S(A^2 R'\bo u^1x^2\o Y^3\o\o \underline{X^2\o d^3\o} T^2)\alpha \\
	&		& \,\,\,\,\, A^3u^3 x^2\t\t Y^3\o\t\t \underline{X^2\t\t d^3\t\t} T^3\t t^2\t h^2D^3\o\\
	&		& \tens a^3\t G^3x^3Y^3\t \underline{X^3}t^3 h^3 D^3\t\\
	&	=	& a^1Y^1 b\o d^1\underline{H^1t^1}D^1\\
	&		& \tens a^2G^1R\bt A^1 R'\bt u^2 x^2\t\o Y^3\o\t\o d^3\o\t\o X^2\t\o \underline{H^2\t\t\o T^3\o t^2\o} h^1D^2\\
	&		& \tens a^3\o G^2R\bo x^1Y^2 b\t d^2X^1\underline{H^2\o T^1}\beta S(A^2 R'\bo u^1x^2\o Y^3\o\o d^3\o\o X^2\o \underline{H^2\t\o T^2})\alpha \\
	&		& \,\,\,\,\, A^3u^3 x^2\t\t Y^3\o\t\t d^3\o\t\t X^2\t\t \underline{H^2\t\t\t T^3\t t^2\t} h^2D^3\o\\
	&		& \tens a^3\t G^3x^3Y^3\t d^3\t X^3\underline{H^3 t^3}h^3 D^3\t\\
	&	=	& a^1Y^1 b\o d^1D^1\\
	&		& \tens a^2G^1R\bt \underline{A^1 R'\bt u^2 x^2\t\o Y^3\o\t\o d^3\o\t\o X^2\t\o} T^3\o h^1D^2\\
	&		& \tens a^3\o G^2R\bo x^1Y^2 b\t d^2X^1T^1\beta S(\underline{A^2 R'\bo u^1x^2\o Y^3\o\o d^3\o\o X^2\o} T^2)\alpha \\
	&		& \,\,\,\,\, \underline{A^3u^3 x^2\t\t Y^3\o\t\t d^3\o\t\t X^2\t\t} T^3\t h^2D^3\o\\
	&		& \tens a^3\t G^3x^3Y^3\t d^3\t X^3 h^3D^3\t\\
	&	=	& a^1Y^1 b\o d^1D^1\\
	&		& \tens a^2G^1R\bt x^2Y^3\o d^3\o X^2\underline{A^1 R'\bt u^2 T^3\o} h^1D^2\\
	&		& \tens a^3\o G^2R\bo x^1Y^2 b\t d^2X^1\underline{T^1\beta S(A^2 R'\bo u^1 T^2)\alpha A^3u^3  T^3\t} h^2D^3\o\\
	&		& \tens a^3\t G^3x^3Y^3\t d^3\t X^3 h^3D^3\t\\
	&	=	& a^1Y^1 b\o d^1D^1 \tens a^2G^1R\bt x^2Y^3\o d^3\o X^2R^-\bt h^1D^2\\
	&		& \tens a^3\o G^2R\bo x^1Y^2 b\t d^2X^1R^-\bo h^2D^3\o \tens a^3\t G^3x^3Y^3\t d^3\t X^3 h^3D^3\t\\
	&	=	& x^1Y^1 b\o y^1X^1  \tens x^2T^1R\bt w^2Y^3\o y^3\o W^2R^-\bt t^1X^2\\
	&		& \tens x^3\o T^2R\bo w^1Y^2 b\t y^2W^1R^-\bo t^2X^3\o \tens x^3\t T^3w^3Y^3\t y^3\t W^3t^3X^3\t\\
	&	=	& \Delta(b\tens 1)
\end{eqnarray*}
\endproof

\begin{example} Recall the structure of $\underline{D^{\phi}(G)}$ from section 3. We can now compute the structure of $\underline{D^{\phi}(G)}\bos D^{\phi}(G)$.

\begin{eqnarray*}
\lefteqn{((g\tens\delta_s)\tens(g'\tens\delta_{s'}))((h\tens\delta_t)\tens(h'\tens\delta_{t'}))=}\\
	&		&	(gg'hg'^{-1}\tens\delta_{gg'tg'^{-1}g^{-1}})\tens(g'h'\tens\delta_{g't'g'^{-1}}) \, \delta_{s,gg'tg'^{-1}g^{-1}}\delta_{s',g'th^{-1}t^{-1}ht'g'^{-1}}\\
	&		&	\theta_{g't'g'^{-1}}(g',h')\theta_{gg'tg'^{-1}g^{-1}}(g,g'hg'^{-1})\gamma_{g'}(g'th^{-1}t^{-1}hg'^{-1},g't'g'^{-1})\\
	&		&	\gamma_{g'}(g'tg'^{-1},g'h^{-1}t^{-1}hg'^{-1})\theta^{-1}_{g'h^{-1}t^{-1}hg'^{-1}}(g',g'^{-1})\\
	&		&	\gamma^{-1}_{g'}(g'h^{-1}t^{-1}hg'^{-1},g'h^{-1}thg'^{-1})\theta_{g'tg'^{-1}}(g',h)\theta_{g'tg'^{-1}}(g'h,g'^{-1})\\
	&		&	\phi(gg'tg'^{-1}g^{-1},g't^{-1}g'^{-1},g'h^{-1}t^{-1}hg'^{-1})\phi^{-1}(g't^{-1}g'^{-1},g'tg'^{-1},g'h^{-1}t^{-1}hg'^{-1})\\
	&		&	\phi^{-1}(gg'tg'^{-1}g^{-1}g't^{-1}g'^{-1},g'th^{-1}t^{-1}hg'^{-1},g't'g'^{-1})
\end{eqnarray*}

\[ \eta(1) = \sum_{s,t\in G} (e\tens\delta_s)\tens(e\tens\delta_t) \, \phi(s^{-1},s,s^{-1}) \]

\begin{eqnarray*}
\lefteqn{\Delta((g\tens\delta_s)\tens(g'\tens\delta_{s'})) =}\\
	&		& \sum_{jk=s}\sum_{ab=t} (kgk^{-1}\tens\delta_j)\tens(kg^{-1}k^{-1}gg'\tens\delta_{kg^{-1}k^{-1}gag^{-1}kgk^{-1}})\tens(g\tens\delta_k)\tens(g'\tens\delta_b)\\
	&		& \gamma_{g'}(a,b)\gamma_g(j,k)\theta^{-1}_j(kgk^{-1},kg^{-1}k^{-1}g)\theta_{kg^{-1}k^{-1}gag^{-1}kgk^{-1}}(kg^{-1}k^{-1}g,g')\\
	&		&	\phi(s,g^{-1}s^{-1}g,g^{-1}sg)\phi^{-1}(j,kg^{-1}k^{-1}j^{-1}kgk^{-1},kg^{-1}k^{-1}jkgk^{-1})\phi^{-1}(k,g^{-1}k^{-1}g,g^{-1}kg)\\
	&		& \phi^{-1}(jkg^{-1}k^{-1}j^{-1}kgk^{-1},kg^{-1}k^{-1}g,g^{-1}jkg)\phi^{-1}(kg^{-1}k^{-1}jkgk^{-1},kg^{-1}k^{-1}g,g^{-1}kg)\\
	&		& \phi(kg^{-1}k^{-1}g,g^{-1}jg,g^{-1}kg)\phi(jkg^{-1}k^{-1}j^{-1}kgk^{-1},kg^{-1}k^{-1}jkgk^{-1},k)\\
	&		& \phi(jkg^{-1}k^{-1}j^{-1}kgk^{-1},kg^{-1}k^{-1}g,ab)\phi(kg^{-1}k^{-1}gag^{-1}kgk^{-1},kg^{-1}k^{-1}g,b)\\
	&		& \phi^{-1}(kg^{-1}k^{-1}g,a,b)\phi^{-1}(jkg^{-1}k^{-1}j^{-1}kgk^{-1},kg^{-1}k^{-1}gag^{-1}kgk^{-1},kg^{-1}k^{-1}gb)
\end{eqnarray*}

\[ \varepsilon((g\tens\delta_s)\tens(g'\tens\delta_{s'})) = \delta_{s,e}\delta_{s',e} \]

\begin{eqnarray*}
\lefteqn{S((g\tens\delta_s)\tens(g'\tens\delta_{s'})) =}\\
	&		&	(g'^{-1}g^{-1}sg^{-1}s^{-1}gg'\tens\delta_{g'^{-1}g^{-1}s^{-1}gg'})\tens(g'^{-1}g^{-1}sgs^{-1}\tens\delta_{g'^{-1}g^{-1}sgs^{-1}s'^{-1}g'})\\
	&		&	\theta^{-1}_{s^{-1}}(g,g^{-1})\gamma^{-1}_g(s,s^{-1})\theta^{-1}_{sg^{-1}s^{-1}gs'^{-1}g^{-1}sgs^{-1}}(sg^{-1}s^{-1}gg',g'^{-1}g^{-1}sgs^{-1})\\
	&		& \gamma^{-1}_{sg^{-1}s^{-1}gg'}(sg^{-1}s^{-1}gs'g^{-1}sgs^{-1},sg^{-1}s^{-1}gs'^{-1}g^{-1}sgs^{-1})\theta_{sg^{-1}s^{-1}gs^{-1}}(sg^{-1}s^{-1}g,g^{-1})\\
	&		& \theta_{sg^{-1}s^{-1}gs'g^{-1}sgs^{-1}}(sg^{-1}s^{-1}g,g')\gamma_{g'^{-1}g^{-1}sgs^{-1}}(g'^{-1}sg^{-1}s^{-1}gg',g'^{-1}g^{-1}sgs^{-1}s'^{-1}g')\\
	&		& \gamma_{g'^{-1}g^{-1}sgs^{-1}}(g'^{-1}g^{-1}s^{-1}gg',g'^{-1}g^{-1}sgsg^{-1}s^{-1}gg')\\
	&		& \theta^{-1}_{g'^{-1}g^{-1}sgs^{-1}g^{-1}s^{-1}gg'}(g'^{-1}g^{-1}sgs^{-1},sg^{-1}s^{-1}gg')\\
	&		& \gamma^{-1}_{g'^{-1}g^{-1}sgs^{-1}}(g'^{-1}g^{-1}sgsg^{-1}s^{-1}gg',g'^{-1}g^{-1}sgs^{-1}g^{-1}s^{-1}gg')\\
	&		& \theta_{g'^{-1}g^{-1}s^{-1}gg'}(g'^{-1}g^{-1}sgs^{-1},sg^{-1}s^{-1})\phi(sg^{-1}s^{-1}gs'g^{-1}sgs^{-1},sg^{-1}s^{-1}g,s'^{-1}g^{-1}sgs^{-1})\\
	&		& \phi(sg^{-1}s^{-1}gs^{-1},sg^{-1}s^{-1}g,g^{-1}sg)\phi(s,g^{-1}s^{-1}g,g^{-1}sg)\phi^{-1}(sg^{-1}s^{-1}g,g^{-1}s^{-1}g,g^{-1}sg)
\end{eqnarray*}

\[ \alpha = \sum_{s,t\in G}(e\tens\delta_s)\tens(e\tens\delta_t) \phi(s^{-1},s,s^{-1}) \]

\[ \beta = \sum_{s,t\in G}(e\tens\delta_s)\tens(e\tens\delta_t) \phi(s^{-1},s,s^{-1})\phi(t^{-1},t,t^{-1}) \]

\begin{eqnarray*}
\phi	& =	& \sum_{g,h,k\in G}\sum_{u,v,w\in G} ((e\tens\delta_g)\tens(e\tens\delta_u))\tens((e\tens\delta_h)\tens(e\tens\delta_v))\tens((e\tens\delta_k)\tens(e\tens\delta_w))\\
			&		&  \,\,\,\,\,\,\, \phi(u,v,w)\phi(g^{-1},g,g^{-1})\phi(h^{-1},h,h^{-1})\phi(k^{-1},k,k^{-1})
\end{eqnarray*}

The isomorphism for this example is

\begin{eqnarray*}
\chi((g\tens\delta_s)\tens(h\tens\delta_t))	&	=	& (gh\tens \delta_s)\tens (h\tens \delta_{g^{-1}s^{-1}gt})\theta_s(g,h)\gamma_h(g^{-1}sg,g^{-1}s^{-1}gt)\\
	&		& \,\,\,  \phi(s,g^{-1}s^{-1}g,g^{-1}sg)\phi^{-1}(sg^{-1}s^{-1}g,g^{-1}sg,g^{-1}s^{-1}gt)
\end{eqnarray*}
\end{example}

\end{document}